\documentclass[a4paper, 11pt]{article}
\usepackage{graphicx,comment,subcaption}
\usepackage{amsfonts}
\usepackage{amsmath}
\usepackage{amssymb}
\usepackage[round]{natbib}
\usepackage{amsthm}
\usepackage{algorithmic}
\usepackage{color,soul}
\usepackage{tikz,url} 
\usepackage{natbib}
\usepackage[colorlinks=true, allcolors=blue]{hyperref}
\usepackage{booktabs}
\usepackage{adjustbox}
\usepackage[ruled,vlined]{algorithm2e}
\usepackage{pgfplots}
 
\usetikzlibrary{positioning}

\newtheorem*{example*}{Example}

\newcommand{\R}{\mathbb R}
\newcommand{\D}{\mathcal D}

\allowdisplaybreaks[1]

\title{Leveraging Contextual Information for Robustness in Vehicle Routing Problems}

\author{
  Ali İrfan Mahmutoğulları\\
  \texttt{irfan.mahmutogullari@kuleuven.be}\\
  KU Leuven, Leuven, Belgium 
  \and
  Tias Guns\\
  \texttt{tias.guns@kuleuven.be} \\
  KU Leuven, Leuven, Belgium 
}

\begin{document}

\maketitle
\begin{abstract}
   We investigate the benefit of using contextual information in data-driven demand predictions to solve the robust capacitated vehicle routing problem with time windows. Instead of estimating the demand distribution or its mean, we introduce contextual machine learning models that predict demand quantiles even when the number of historical observations for some or all customers is limited. We investigate the use of such predicted quantiles to make routing decisions, comparing deterministic with robust optimization models. Furthermore, we evaluate the efficiency and robustness of the decisions obtained, both using exact or heuristic methods to solve the optimization models. Our extensive computational experiments show that using a robust optimization model and predicting multiple quantiles is promising when substantial historical data is available. In scenarios with a limited demand history, using a deterministic model with just a single quantile exhibits greater potential. Interestingly, our results also indicate that the use of appropriate quantile demand values within a deterministic model results in solutions with robustness levels comparable to those of robust models. This is important because, in most applications, practitioners use deterministic models as the industry standard, even in an uncertain environment. Furthermore, as they present fewer computational challenges and require only a single demand value prediction, deterministic models paired with an appropriate machine learning model hold the potential for robust decision-making.
\end{abstract}


\section{Introduction}
The increasing availability of data in various industries has had a profound impact on the applications of operations research (OR) in recent years. Examples, as mentioned in \cite{bertsimas_data-driven_2018}, include transaction data stored by retailers, order trends observed by suppliers, and real-time power usage data in energy markets.  This vast amount of transaction data is being archived by retailers, amounting to terabytes of information. Similarly, suppliers continuously monitor and track order patterns throughout their supply chains. Furthermore, energy markets have access to data, including global weather information, historical demand profiles and real-time power consumption data. This trend has sparked research focused on integrating machine learning (ML) and OR, with a strong emphasis on data analysis.

On the other hand, in this rapidly evolving and data-rich environment, practitioners face the challenge of ensuring the feasibility and quality of their decisions. Robust optimization models provide a viable solution by generating tractable problems and delivering solutions that are comparable to those obtained through alternative methods, given the appropriate uncertainty sets for the uncertain problem parameters (see, for example, \cite{ben-tal_robust_2002} and \cite{bertsimas_theory_2011} for a deeper discussion of robust optimization models). These models prove particularly valuable when dealing with problems where no distributional information is available for the problem parameters.

Another recent development in the OR community is the use of contextual information in decision-making. Contextual information, also known as features or covariates, refers to auxiliary data that exhibits a correlation with the actual problem parameters. For instance, the demand for a specific product by a customer is influenced by the customer's profile. Thus, in practical applications, the price of a new product is determined by observing the customers' features and predicting their demands rather than observing the true demand of the customers. Another example is a real-life scenario encountered by a large logistics service provider (LSP), which also motivates our study on the utilization of contextual information for addressing the robust capacitated vehicle routing problem with time windows (RCVRPTW). The LSP provides routing services to over 6,000 customers and the demand values for these customers exhibit uncertainty. Notably, while more than half of all customers have only one recorded demand observation, a mere 3\% of customers possess over 30 demand observations. This contrast renders accurate forecasting of each customer's demand value impossible. However, the LSP possesses contextual information pertaining to customer orders, including origin, destination and packaging type. It is known that contextual information exhibits a correlation with demand values. Occasionally, the LSP encounters new customers lacking a demand history; nonetheless, certain contextual information about these new customers remains available, such as the location of the customers, the types of goods they ship and possibly some contextual information about the orders, such as packaging types. Consequently, the LSP aims to devise optimal vehicle routes by effectively predicting customer demand, thereby enabling the provision of services at a minimal cost. Moreover, the LSP would like to ensure the robustness of vehicle routes to avoid infeasible solutions after realizing the actual value of customer demands since they require expensive recourse actions.
 
In order to make routing decisions for its vehicles, the LSP uses demand estimates of customers, denoted by the vector $\tilde{d}$. However, due to the uncertainty of demand values, some of the routes obtained with $\tilde{d}$ may become infeasible after the actual demand realization $d$. In such cases, the vehicles have to diverge from their routes at some point for replenishment or unloading before serving the remaining customers on the route. Consider a given route $(0,i_1,i_2,\dots,i_K,0)$ where 0 denotes the depot and $i_1,i_2,\dots,i_K$ are the customers served on this route. If the traveling cost from $i$ to $j$ is $c_{ij}$ then the initial cost of the route is $c_{0i_1} + c_{i_1i_2} + \cdots + c_{i_{K-1}i_{K}} + c_{i_K0}$. Now, assume that for a specific demand vector $d$, a vehicle with capacity $Q$ can serve the first $k < K$ customers on the route but is unable to serve the next customer; that is, $d_1 + \cdots + d_k \leq Q$ but $d_1 + \cdots + d_{k+1} > Q$. In this scenario, the vehicle should visit the depot instead of customer $i_{k+1}$ immediately after serving customer $i_k$ as a recourse action, incurring an additional cost of $c_{i_{k}0} + c_{0i_{k+1}} - c_{i_{k} i_{k+1}}$ due to this detour. Note that it is also possible to have multiple detours on a single route. Then, the overall operational cost is the sum of the initial route cost and all the required recourse costs due to detours. Another important measure for evaluating the robustness of solutions is the amount of time window violations and the percentage of customers that cannot be served within their specified time windows. Although the planned routes may satisfy the service time window constraints of all customers, delays caused by detours can lead to potential violations for some of the customers and result in additional costs due to customer dissatisfaction. 

This fact raises the question of how to obtain demand predictions, $\tilde{d}$, that lead to the best decisions in terms of total operating cost and robustness of the solution against time window violations. One approach is to estimate the demand distribution for each customer separately. However, this may not always be possible due to limited demand realization data for each customer, as in the LSP case. Alternatively, one can make the routing decisions based on the mean or average demand predictions for the customers. However, particularly when demand variability is high, relying solely on predicted means can result in routes with a large number of customers. This can be problematic, especially if the actual demand realizations significantly exceed the mean demand predictions, as it may lead to infeasible routes and hence expensive recourse actions. In this case, predicting an overestimation, e.g., a value above the expected demand, will create a buffer for demand realizations where the true value is higher than the expected value. A technical way to achieve this is to predict not the expected value but to do a quantile regression at a specific quantile value. For instance, for a 0.6 quantile value, the prediction will be such that, in expectation, 60\% of the demand realizations will be below the predicted value. Different quantile values will hence create different amounts of buffer in its prediction while automatically taking the variance observed in the data into account without any assumption of the distribution (e.g., the amount of buffer at 0.6 quantile for a highly volatile customer will be much higher than that of a low-variance customer for the same expected mean).

%A promising approach is the use of demand quantiles, which enables the decision-maker to fine-tune the prediction level. For instance, by selecting the 0.6 quantile, which represents values greater than 60\% of all demand realizations, the decision-maker can create a buffer, reducing the need for recourse actions that may lead to additional costs and potential time window violations. 
%provides a buffer for actual demand realizations, enhancing the solution's robustness.

Therefore, in this paper, we propose the use of predicted demand quantiles obtained from the contextual information of the customers to make robust decisions in vehicle routing applications. The use of quantiles offers two distinct advantages. Firstly, it does not require any assumptions regarding the distribution of customer demand; instead, it directly predicts the quantile values. Secondly, in contrast to mean predictions, it enables the decision-maker to craft solutions with a degree of robustness by predicting values beyond the mean or median, such as a quantile for the demands. In the deterministic formulation of the problem, a single demand prediction is made for each customer. However, in the robust formulation where demand uncertainty is explicitly represented in the model, multiple predictions are required (e.g., two values are needed to define an interval). These predictions are utilized to create the uncertainty set in the robust formulation, which includes all possible realizations of the demand values.

Therefore, in this paper, we explore the following research questions:

\begin{itemize}
\item Can quantile predictions for customer demands be used in the context of deterministic or robust capacitated vehicle routing problems with time windows (CVRPTW) when the number of past observations is limited for some or all customers?

\item Can the utilization of contextual information from existing and new customers facilitate the prediction of uncertainty sets for demands in  RCVRPTW, thereby enabling the formulation of robust decisions?

\item  Is it possible for non-expert practitioners to obtain robust solutions by solving the deterministic CVRPTW using appropriately selected quantile predictions?
\end{itemize}
 
We provide answers to these questions and contribute to the existing literature in the following ways:

\begin{itemize}
    \item[(i)] We propose contextual prediction models that enable the estimation of customer demand quantiles or means, even in cases where the number of past observations is limited for some or all customers. To the best of our knowledge, this study is the first to consider contextual prediction models in vehicle routing problems. These models are designed to make predictions for both existing and new customers based on contextual information such as the type of good, location of the customer, etc. We explore both linear and non-linear prediction models, where the non-linear models are defined by neural networks. The predicted quantile values facilitate the solution of CVRPTW and its robust counterpart, RCVRPTW, enhancing the decision-making process in terms of cost and robustness.  

    \item[(ii)]  We assess the robustness of decisions obtained using mathematical models or heuristics when quantile demand predictions are used. For small problem instances, mathematical models can be effectively solved using commercial mixed-integer programming (MIP) solvers. However, for larger instances, solvers are not able to find the optimal solution at reasonable solution times, and therefore heuristic methods can be employed to solve the problem efficiently. Therefore, we extend the standard implementation of an adaptive large neighborhood search (ALNS) based heuristic for the deterministic problem CVRPTW to its robust counterpart RCVRPTW. 

    \item[(iii)] We demonstrate that using appropriate quantile demand values in a deterministic formulation yields solutions that are almost at least as robust as those obtained from the robust models in the data-scarce setting. This finding carries significance since the formulation and solution methods for the robust model may pose challenges for practitioners, while the deterministic formulation may be preferred by them. Also, the deterministic models are industry standards, as they are widely utilized in practice. 
\end{itemize}

The rest of the paper is structured as follows: In Section \ref{sec:literature}, we provide a comprehensive review of the related literature. In Section \ref{sec:models}, we present the notation, formal definition of the problem and mathematical models for both the CVRPTW and its robust counterpart, the RCVRPTW. In Section \ref{sec:alns}, we present an ALNS heuristic for solving large instances that are not solvable in a reasonable time with the mathematical models presented in Section \ref{sec:models}. Prediction of quantiles using contextual information is explained in Section \ref{sec:prediction}. The results of our computational experiments are discussed in Section \ref{sec:experiments}. Finally, we present our concluding remarks and possible extensions in Section \ref{sec:conclusion}.

\section{Literature Review}\label{sec:literature}
The literature review section is divided into four subsections to provide a comprehensive overview of the relevant research. The first subsection focuses on the related literature concerning robust vehicle routing problems. The second subsection covers studies that explore the generation of uncertainty sets in the context of robust optimization problems. The third subsection is dedicated to studies that specifically explore the utilization of contextual information in decision-making processes. Lastly, the fourth subsection presents studies that focus on quantile prediction methods.

\subsection{Robust Vehicle Routing Problem}
Due to the extensive body of literature on vehicle routing problems (VRP), numerous studies have been dedicated to this particular OR application. However, for the purpose of this literature review, our focus will primarily be on the robust VRP literature as it pertains specifically to the problem addressed in this paper. For a detailed discussion on general VRPs, we refer interested readers to \cite{laporte_vehicle_1992} and \cite{toth_vehicle_2002}. 

While uncertainty in the context of VRPs has also been studied in several works, the majority of the existing work focuses on stochastic programming models. These models typically assume that the uncertain problem parameters follow known probability distributions. However, this assumption may not hold true for many real-life applications. \cite{oyola_stochastic_2018} present detailed models and discussions on stochastic VRPs. While stochastic VRP models have received significant attention, robust optimization approaches in the VRP context have been less explored.

To the best of our knowledge, \cite{sungur_robust_2008} are the first to address robust VRPs with uncertain demand values. They propose mathematical models with Miller-Tucker-Zemlin (MTZ) constraints that incorporate different types of uncertainty sets, namely convex hull, box, and ellipsoidal. They show that solving these robust models is computationally equivalent to solving the deterministic version of the problem with augmented demand values based on the shape of the uncertainty set. \cite{sungur_robust_2008} focus on worst-case robustness, meaning that the solutions obtained from their models are feasible for all demand realizations within the defined uncertainty set. This overly-conservative approach ensures that the solutions are always feasible under any realization of the demand values. In a more general setting, \cite{bertsimas_robust_2003} prove that any robust combinatorial optimization problem can be solved by solving  its deterministic counterpart a polynomial number of times if the uncertain parameters are in the objective. Thus, robust VRP with uncertain arc cost values can be solved by solving $n+1$ deterministic VRPs where $n$ is the number of arcs. They consider a budgeted uncertainty  set where a parameter $\Gamma$ adjusts the degree of robustness of the solution. 

\cite{gounaris_robust_2013} developed exact robust counterparts of the deterministic capacitated VRP (CVRP) using different formulations, including two-index, Miller-Tucker-Zemlin, commodity flow, and vehicle assignment formulations. For each formulation, the authors established necessary and sufficient conditions to define a robust feasible set of routes. Later, \cite{gounaris_adaptive_2016} propose a meta-heuristic for solutions of robust CVRP under polyhedral demand uncertainty. The key focus of their meta-heuristic algorithm is to ensure the robust feasibility of the generated route  to ensure that they remain feasible for all demand realizations.

\cite{munari_robust_2019} address RCVRPTW where both demand and travel time values are uncertain. Firstly, they introduce a compact mathematical formulation for the RCVRPTW, leveraging the recursive dynamic programming representation of budgeted uncertainty sets. This formulation offers a simpler representation of robustness 
than other existing formulations based on dualization of the uncertainty sets. Moreover, it demonstrates improved overall performance when solved using commercial solvers. \cite{munari_robust_2019}  also propose a branch-and-price method based on the set partition formulation of the problem. The results of the computational experiments show the importance of making robust decisions in the existence of uncertainty for real-life applications.

Finally, we refer an interested reader to \cite{hasenbein_robust_2010} for a tutorial paper on robust VRP models.

\subsection{Generating Uncertainty Sets in Robust Problems}
Traditionally, uncertainty sets in robust optimization problems are assumed to be given and available to the decision maker. However, recent studies also showed that it is also possible to construct these sets based on the available data the decision maker has. 

\cite{goldfarb_robust_2003} propose generating uncertainty sets using market data and a parameter that controls the size of the uncertainty set in the context of a portfolio optimization problem. They assume that the returns of random assets are linearly dependent on the factors that drive the market. In a related study, \cite{tulabandhula_robust_2014} adopt a similar approach by focusing on designing uncertainty sets for a broad class of decision-making problems while making minimal assumptions about the distribution. Their work introduces uncertainty set design based on statistical learning theory. Similar to our work, \cite{tulabandhula_robust_2014} propose a conditional quantile prediction method to generate uncertainty sets for robust optimization problems using contextual information. By utilizing conditional quantiles, they provide a robustness guarantee within the robust formulation. However, they do not explicitly mention specific applications or provide computational results for this approach to generating uncertainty sets.

In their work, \cite{bertsimas_data-driven_2018} introduce a novel framework for constructing uncertainty sets in robust optimization using data and hypothesis tests. Their approach leads to uncertainty sets that offer a probabilistic guarantee and are typically smaller than sets derived from limited data. This results in less conservative models compared to traditional robust approaches, while preserving the same level of robustness guarantees.  

\cite{ohmori_predictive_2021} presents an algorithm that uses minimum volume ellipsoids enclosing $k$-nearest neighbors in feature space as uncertainty sets. The algorithm employs a nonparametric prediction model, allowing for flexibility in handling diverse data patterns without making strong distributional assumptions.  This approach provides robustness against prediction errors. Moreover, the use of ellipsoidal uncertainty sets benefits from extensive research on efficient algorithms in robust optimization.

In their recent study, \cite{goerigk_data-driven_2023} applied unsupervised deep learning methods to generate uncertainty sets for robust optimization. They utilized machine learning techniques to uncover anomalies within the data, leading to the creation of non-convex uncertainty sets. This approach allows for the derivation of robust solutions by formulating the adversarial problem as a convex quadratic mixed-integer program.  In contrast, we employ contextual information about the customers to predict multiple quantiles, which are used for defining the uncertainty sets in the RCVRPTW model.

\subsection{Contextual Optimization} 

The emerging paradigm of contextual optimization has attracted significant interest from both the ML and OR communities. In contextual optimization problems, decision-making is done without directly observing the uncertain problem parameters, such as objective function or constraint coefficients. Instead, decisions are made based on available contextual information, also referred to as features or covariates.

In their respective works, \cite{ban_big_2018} and \cite{zhang_optimal_2023} tackle contextual newsvendor problems with different approaches. \cite{ban_big_2018} introduces algorithms based on Empirical Risk Minimization and Kernel-weights Optimization, enabling decision-making solely based on observed features and eliminating the need for a separate demand estimation step. They also provide performance bounds specifically for the feature-driven decision case. On the other hand, \cite{zhang_optimal_2023} proposes a distributionally robust framework utilizing the Wasserstein distance to handle uncertainties associated with unseen data. Their focus is on preventing ill-defined policies when faced with data that deviates from the observed historical distribution. 

In the study conducted by \cite{ban_dynamic_2019}, the authors introduce the residual tree method as a means to forecast demand for new products using contextual information. By incorporating covariates and generating multiple demand values within a scenario tree structure, their approach improves the accuracy of demand forecasting and allows decision-makers to assess the robustness of their decisions under various scenarios. Similarly, in the research by \cite{kannan_data-driven_2022}, the authors predict point estimates and out-of-sample residuals of leave-one-out models for generating scenarios, enabling the use of sample average approximation (SAA) to approximate the true distribution of problem parameters.   

In their study, \cite{bertsimas_predictive_2020} explore the use of contextual information in prediction and decision-making. They propose a methodology that incorporates predictive machine learning methods applied to feature values, allowing decisions to be made based on a weighted combination of historical values. The authors consider different approaches, such as $k$-nearest-neighbors, kernel methods, local linear models, prediction trees, and ensembles, to determine these weights. Furthermore, the paper investigates decision-dependent parameters and examines the asymptotic behavior of the prescriptions.  Similarly, in the work of \cite{lin_data-driven_2022}, feature information of existing products is leveraged to optimize the order quantity of a new product based on its similarity to the existing products. By utilizing contextual information, both studies aim to improve the decision-making processes and achieve better performance in their respective domains.

In contrast to the studies that separate prediction and optimization, decision-focused learning (DFL) approaches integrate prediction and optimization into a single pipeline. DFL distinguishes itself from other contextual methods by incorporating the decision error, which is dependent on the downstream optimization problem, into the training process. One example is the work by \cite{bertsimas_optimal_2019}, where they introduce the concept of optimal prescriptive trees where the loss to be minimized is a convex combination of prediction and prescription errors. By obtaining prescriptions based on feature values, this approach offers benefits such as interpretability, scalability, generalizability and competitive performance when compared to alternative methods. 

In the study conducted by \cite{chung_decision-aware_2022}, contextual optimization is applied to optimize a healthcare supply chain. The authors propose a method that learns a weighting of the points in the training data to minimize decision loss, effectively incorporating contextual information into the decision-making process. 

In this work, we utilize contextual information about the customers to enhance robust decision-making in the vehicle routing setting. To the best of our knowledge, no existing literature has explored contextual optimization for routing problems. An interested reader may refer to \cite{sadana_survey_2023} for a general discussion on contextual optimization models. 

%In their work, \cite{elmachtoub_smart_2022} introduce the smart ``predict, then optimize" (SPO) framework, which focuses on minimizing the decision loss directly in problems where uncertain parameters appear linearly in the objective. To address the challenge of vanishing gradients, they propose a convex surrogate loss function called SPO+ that is statistically consistent with SPO. The authors also demonstrate the promising performance of SPO in cases of model misspecification. Furthermore, the literature explores various extensions and improvements, including the use of contrastive losses \cite{mulamba_contrastive_2021}, learning-to-rank techniques \cite{mandi_decision-focused_2022}, smoothing the objective with a quadratic term \cite{wilder_melding_2019} , and applying interior point methods with a log barrier term \cite{mandi_interior_2020}. In a recent study by \cite{estes_smart_2023}, the use of contextual information is considered in two-stage linear programming problems, further expanding the applications of contextual optimization to two-stage stochastic optimization in decision-making.

%\cite{mandi_decision-focused_2023} survey the state-of-the-art DFL techniques as well as various use cases and possible research directions.  

\subsection{Quantile Regression}
In statistics, the $\beta$ quantile of a random variable $X$ corresponds to the minimum (or infimum) $x$ value such that the cumulative distribution function $F(x)$ is equal to or greater than $\beta$. For instance, the $0.5$ quantile, also known as the median, is commonly used in statistics. It represents the value for which half of the data falls below and half falls above. Quantile regression, introduced by \cite{koenker_regression_1978}, is a statistical technique that focuses on predicting specific quantiles of a random variable using observed covariates. Unlike traditional regression methods that estimate the conditional mean, quantile regression provides estimates for different quantiles of the response variable.

Quantile regression models aim to minimize a tilted absolute value function, often referred to as the check function or pinball loss function. This function places different weights on observations that fall above or below the predicted quantile, which helps capture the asymmetric nature of the quantiles, unlike the classical mean squared error (MSE). In Section \ref{sec:prediction}, we provide a more detailed explanation of quantile regression.

The relationship between covariates and predicted quantiles can be linear (\cite{koenker_quantile_2001}) or nonlinear, e.g., defined by a neural network (\cite{hatalis_novel_2019,brando_deep_2022}). When multiple quantiles are predicted simultaneously, a challenge known as quantile crossing may arise, where the estimated quantile values exhibit an undesirable reversal in their order as the quantile level increases. For instance, \cite{dai_non-crossing_2022} propose using smoothing techniques to eliminate crossing quantiles and ensure a coherent ordering of the predicted quantiles. This phenomenon is not the main focus of our study and we adopt a simple strategy to mitigate it as explained in Section \ref{sec:prediction}. 

\section{Mathematical Model} \label{sec:models}
We will first present the mathematical model for the deterministic CVRPTW for the sake of clarity. Then, we will extend this model to its robust counterpart. The problem involves a depot and a set of customers to serve, denoted by $\mathcal{C} = \{1, \ldots, n\}$. The CVRPTW can be represented using a graph $\mathcal{G}=(\mathcal{N},\mathcal{A})$, where $\mathcal{N} = \mathcal{C} \cup \{0,n+1\}$ represents the set of nodes and $\mathcal{A} = \{(0,j) : j \in \mathcal{C}\} \cup \{(i,j) : i,j \in \mathcal{C}, i \neq j\} \cup \{(i,n+1) : i \in \mathcal{C}\}$ represents the set of arcs. Both nodes 0 and $n+1$ correspond to the depot node.

A vehicle's route must begin at node 0, visit a subset of customers and finally return to node $n+1$. The time and cost required to traverse arc $(i,j)$ are denoted by $t_{ij}$ and $c_{ij}$, respectively. We assume that the cost is proportional to the travel time. Each customer $i$ has a predicted demand denoted by $\tilde{d}_i$, where $\tilde{d}_0 = \tilde{d}_{n+1} = 0$. Recall that the routing decisions are made based on a predicted demand vector $\tilde{d}$, not on the actual customer demands. Moreover, the vehicles have a capacity denoted by $Q$ and the total demand of the customers served on the same route cannot exceed this capacity. To serve a customer $i$, it must be visited within a specific time window $[T^{min}_i, T^{max}_i]$. Additionally, each customer $i$ has a service time denoted by $s_i$.

The binary decision variable $x_{ij}$ indicates if arc $(i,j)$ is used or not. The continuous decision variable $u_i$ indicates the amount of demand served right after visiting node $i$ by the vehicle to which node $i$ is assigned. Here, we formulate the problem as a collection  Finally, the continuous decision variable $w_i$ represents the visiting time of node $i$. If the predicted demands of customers are given by $\{\tilde{d}_i, i \in \mathcal{C}\}$, we can use the following deterministic model for CVRPTW:
\begin{align}
\min_{x,u,w} &\;\;   \sum_{(i,j) \in \mathcal{A}} c_{ij}x_{ij},  \label{eq:det-obj}\\
\text{s.t.} &\;\; \sum_{(i,j) \in \mathcal{A}} x_{ij} = 1, \;\; \forall j \in \mathcal{C} \label{eq:det-visit} \\
&\;\; \sum_{(i,j) \in \mathcal{A}} x_{ij} - \sum_{(j,i) \in \mathcal{A}} x_{ji} = 0, \;\; \forall i \in \mathcal{C}  \label{eq:det-balance}\\
&\;\; \sum_{(0,j) \in \mathcal{A}} x_{ij} - \sum_{(i,n+1) \in \mathcal{A}} x_{i(n+1)} = 0,  \label{eq:det-vehiclenumber}\\
&\;\; u_{j} \geq u_{i}  + \tilde{d}_j x_{ij} -Q(1-x_{ij}),  \;\; \forall (i,j) \in \mathcal{A}  \label{eq:det-MTZ-demand}\\
&\;\; \tilde{d}_{i} \leq u_{i} \leq Q, \;\; \forall i \in \mathcal{C} \label{eq:det-capacity}\\ 
&\;\; w_{j} \geq w_{i} + (s_i + t_{ij})x_{ij }-T(1-x_{ij}),  \;\; \forall (i,j) \in \mathcal{A} \label{eq:det-MTZ-time}\\ 
&\;\;  T^{min}_i \leq w_i \leq T^{max}_i, \;\; \forall i \in \mathcal{C} \label{eq:det-time-window}\\ 
&\;\; x_{ij} \in \{0,1\},  \;\;  \forall (i,j) \in \mathcal{A} \nonumber\\
&\;\; u_i \geq 0, w_i \geq 0, \;\; \forall i \in \mathcal{C}  \nonumber
\end{align}
The objective \eqref{eq:det-obj} minimizes the total cost. Constraints \eqref{eq:det-visit},\eqref{eq:det-balance} and \eqref{eq:det-vehiclenumber} define the routes of vehicles. Constraints \eqref{eq:det-MTZ-demand} and \eqref{eq:det-capacity} ensure vehicle capacities are not exceeded for the predicted demand values. Time window restrictions are given by Constraints \eqref{eq:det-MTZ-time} and \eqref{eq:det-time-window} where $T$ is a large positive number. 
 
For the robust counterpart, we assume that the demand values belong to intervals, that is, the demand of customer $i$ is in interval $[\overline{d}_i - \hat{d}_i, \overline{d}_i + \hat{d}_i]$ for some base value $\overline{d}_i$ and deviation $\hat{d}_i$. We can equivalently write the demand of customer $i$ as $\overline{d}_i + \xi_i\hat{d}_i$ with $\xi_i \in [-1,1]$ for $i \in \mathcal{C}$ and consider a polyhedral uncertainty set \begin{equation} \label{eq:polyhedral_set}
    D =  \{d \in \R^{|\mathcal{C}|}:  d_i = \overline{d}_i + \xi_i\hat{d}_i, \sum_{i \in \mathcal{C}}\xi_i \leq \Gamma, \xi_i \in [0,1], i \in\mathcal{C} \}
\end{equation}

In Equation \eqref{eq:polyhedral_set}, the parameter $\Gamma$ is called the budget of uncertainty. Indeed, the parameter $\Gamma$ determines the level of risk aversion of the decision-maker. A larger value of $\Gamma$ corresponds to a larger budget of uncertainty, indicating that the decision-maker is more risk-averse. We also assume that all demand values are positive and the worst-case scenario occurs when the demand for customer $i \in \mathcal{C}$ is $\overline{\overline{d}}_i = \overline{d}_i + \hat{d}_i$. A feasible route satisfies Constraints \eqref{eq:det-MTZ-demand} and \eqref{eq:det-capacity} for all possible demand realizations within the uncertainty set $D$. Previous work by \cite{munari_robust_2019} employed a recursive dynamic programming approach to represent the robust constraints when the budget of uncertainty $\Gamma$ is an integer. In our study, we adopt a similar methodology following their approach.

Let $(0,\ldots, i, j, \ldots, n+1)$ be a route and $u_{j\gamma}$ be the load of the vehicle after serving customer $j$ if the worst $\gamma$ customers have reached their worst case demand realization before $j$. Then
\begin{equation}\label{eq:robustness} 
u_{j\gamma} = \max\left\{u_{i\gamma} + \overline{d}_j,u_{i(\gamma-1)} + \overline{\overline{d}}_j \right\} \;\; \forall \gamma \in \{1,\ldots,\Gamma\}, j \in \mathcal{C }
\end{equation}
\begin{equation*} 
    u_{j0} = u_{i0} + \overline{d}_i
\end{equation*}
\begin{equation*} 
    u_{0\gamma} = 0 \;\; \forall \gamma \in \{0,\ldots,\Gamma\}
\end{equation*} 
where the former term in the max operator of \eqref{eq:robustness} corresponds to the case when the demand of customer $j$ is $\overline{d}_j$ and the latter one corresponds to the case when it is $\overline{\overline{d}}_j$. Thus, the following model can be used to solve RCVRPTW where 
the equation \eqref{eq:robustness} is represented by a set of linear equations. \begin{align}
\min_{x,u,w} &\;\;   \sum_{(i,j) \in \mathcal{A}} c_{ij}x_{ij},  \label{eq:rob-obj} \\ 
\text{s.t.} &\;\; \eqref{eq:det-visit} - \eqref{eq:det-vehiclenumber}, \eqref{eq:det-MTZ-time}, \eqref{eq:det-time-window} \nonumber \\
&\;\; u_{j\gamma} \geq u_{i\gamma}  + \overline{d}_j x_{ij} -Q(1-x_{ij}),  \;\; \forall (i,j) \in \mathcal{A}, \forall \gamma \in \{0,\ldots,\Gamma\}\label{eq:rob-MTZ-demand} \\
&\;\; u_{j\gamma} \geq u_{i(\gamma-1)}  + \overline{\overline{d}}_j x_{ij} -Q(1-x_{ij}),  \;\; \forall (i,j) \in \mathcal{A},  \forall \gamma \in \{1,\ldots,\Gamma\}\label{eq:rob-MTZ-gamma} \\
&\;\; \overline{d}_{i} \leq u_{i\gamma} \leq Q, \;\; \forall i \in \mathcal{C}, \forall \gamma \in \{0,\ldots,\Gamma\} \label{eq:rob-capacity}\\
&\;\; x_{ij} \in \{0,1\},  \;\;  \forall (i,j) \in \mathcal{A} \nonumber\\
&\;\; u_{i\gamma} \geq 0, \;\; \forall i \in \mathcal{C},  \forall \gamma \in \{0,\ldots,\Gamma\}  \nonumber \\
&\;\; w_i \geq 0, \;\; \forall i \in \mathcal{C}  \nonumber
\end{align}

In the robust formulation, constraints \eqref{eq:rob-MTZ-demand} and \eqref{eq:rob-MTZ-gamma} ensure that equation \eqref{eq:robustness} holds if customer $j$ is visited right after customer $i$. Constraint  \eqref{eq:rob-capacity} ensures that all routes are feasible in a robust sense since $u_{i\gamma} \leq Q, \forall i \in \mathcal{C}, \forall \gamma \in \{0,\ldots,\Gamma\}$ implies that the vehicle capacities are respected for all demand realizations defined by the uncertainty set \eqref{eq:polyhedral_set}. Also, note that the deterministic model requires a demand prediction $\tilde{d}_i$ for each customer $i \in \mathcal{C}$. However, the robust model requires a base demand $\overline{d}_i$ and a worst-case demand value $\overline{\overline{d}}_i$ for each customer.  

\section{Adaptive Large Neighbourhood Search} \label{sec:alns}
Although the mathematical models presented in the previous section are compact models with two-index variables, they are not able to solve large instances in a reasonable time using commercial MIP solvers. \cite{munari_robust_2019} state that these formulations can solve small instances with 25 customers optimally but can only find feasible solutions for problems with 100 customers in running times of one hour. Therefore, in this section, we present a general ALNS framework for the solution of larger instances.

ALNS is first proposed by \cite{ropke_adaptive_2006} as an extension of the large neighbor search (LNS) method of \cite{shaw_using_1998} as an efficient meta-heuristic search algorithm. The LNS heuristic can search large neighborhoods of a given solution 
using a pair of destroy and repair operations. Unlike LNS, a set of destroy and repair operators are available in ALNS and these operators are chosen based on their performances in the previous iterations.  We use the ALNS heuristic given in Algorithm \ref{alg:ALNS} for solving CVRPTW or its robust counterpart RCVRPTW. 

\begin{algorithm}[ht]
    \SetAlgoLined
    \KwIn{A list of customers $\mathcal{C}$, \\ $\quad$ $\quad$ $\;$ $\;$ Travel cost $c_{ij}$ and time $t_{ij}$ for each $i,j \in \mathcal{C} \cup \{0\}$,  \\  $\quad$ $\quad$$\quad$
    Vehicle capacity $Q$, \\  $\quad$ $\quad$$\quad$  
    Service time windows $[T^{min}_{i},T^{max}_i]$ for each customer $i \in \mathcal{C}$, \\  $\quad$ $\quad$$\quad$  
    A predicted demand vector $\tilde{d}$ for CVRPTW, \\  $\quad$ $\quad$$\quad$  
    (or a base demand vector $\overline{d}$ and a worst-case demand vector $\overline{\overline{d}}$ for RCVRPTW,) \\  $\quad$ $\quad$$\quad$ 
    Set of destroy operators $\Psi^-$, \\  $\quad$ $\quad$$\quad$  
    Set of repair operators $\Psi^+$, \\  $\quad$ $\quad$$\quad$ 
    Initial operator weight vector $W$.}
 
    create an initial feasible solution $s_{\text{min}} = s$\;
    
    \While{stopping criteria not met}{
       
            Select $rep \in \Psi^+$, $des \in \Psi^-$ according to probabilities $p$ obtained from $W$\;
            $s' \leftarrow rep(des(s))$\;
            
            \If{accept($s$, $s'$)}{
                $s \leftarrow s'$\;
                \If{$c(s) < c(s_{\text{min}})$}{
                $s_{\text{min}} \leftarrow s$\;
            }
            }

        Update the weight vector $W$\;
    }
    
    \Return{$s_{\text{min}}$}\;
    
    \caption{A general ALNS framework for CVRPTW and RCVRPTW.} \label{alg:ALNS}
\end{algorithm}

The algorithm starts with a set of routes obtained by a greedy algorithm where the routes are generated by adding the nearest unassigned customer to a route as long as the vehicle capacity is not exceeded. If it is not possible to serve the next customer with the current vehicle capacity, a new route starts from the depot. At each iteration of ALNS with the set of destroy and repair operators $\Psi^-$ and $\Psi^+$, respectively, a destroy operator $des \in \Psi^-$ and a repair operator $rep \in \Psi^+$ are selected to create another feasible solution $s'$ to the problem given the existing solution $s$. To improve the exploration, even if the new solution $s'$ is worse than $s$ in terms of the objective, we can make an update $s \leftarrow s'$ using the  $accept$ function. The $accept$ function in our heuristic is implemented using a simulated annealing mechanism, facilitating enhanced solution exploration at the beginning of the optimization process; that is, even if $s'$ is worse than $s$ in terms of the objective value, it is still accepted with probability $e^{\frac{o(s) - o(s')}{temp}}$ where $temp$ is the temperature parameter and $o(s)$ is the objective function value of solution $s$. The temperature parameter is decreased at each iteration to simulate cooling and reduce exploration in the later iterations of the algorithm. ALNS is stopped after running for a fixed amount of time and the best solution found in terms of initial cost is recorded.  

We use two destroy operators, $\Psi^- = \{random\_removal, string\_removal\}$. The $random\_removal$ operator removes a predefined number of customers from the routes of a feasible solution randomly. The $string\_removal$ operator removes partial routes around a randomly selected customer. As for repair operators, we utilize $\Psi^+ = \{greedy\_repair, regret\_repair\}$. The $greedy\_repair$ approach inserts unassigned customers into existing routes in a greedy manner. The $regret\_repair$ operator first creates a list of unassigned customers according to their decreasing regrets. Regret is determined by calculating the cost difference between assigning a customer to its best position and its second-best position. The $regret\_repair$ operator then inserts these unassigned customers from the list into the solution in a greedy manner, following the order of the list.

ALNS has demonstrated its effectiveness in addressing various routing-type problems, as shown by recent studies such as \cite{hemmelmayr_adaptive_2012, adulyasak_optimization-based_2012, demir_adaptive_2012}. Moreover, several open-source libraries have been developed for the implementation of ALNS. In our study, we use the \texttt{ALNS} library, which is a Python implementation of the ALNS heuristic, provided by \cite{Wouda_Lan_ALNS_2023}.

In order to obtain solutions for both mathematical models presented in Section \ref{sec:models} and the heuristic introduced in this section, some predictions of customer demand values are needed.  In the next section, we will present contextual demand prediction methods with a special emphasis on quantile predictions.

\section{Contextual Quantile Prediction} \label{sec:prediction} 
Assume that we have the demand history $\D_i = \{d_{i1}, d_{i2}, \ldots, d_{in_{i}}\}$ and an assumption for the demand distribution for each customer $i \in \mathcal{C}$. In that case, we can predict the demand distribution of each customer individually when the number of observations $n_i$ is large. For example, the value $\overline{\mu}_i = \frac{1}{n_i} \sum_{l=1}^{n_i}d_{il}$ can be used to predict the mean, and $\sum_{l = 1}^{n_i} \frac{1}{n_i-1}(d_{il} - \overline{\mu}_i)^2$ can be used to predict the variance for the demand distribution of customer $i$. In this case, it becomes possible to sample demand values from the predicted distribution and use them in the decision-making process, such as sample average approximation for a stochastic optimization setting (see, for example, \cite{kleywegt_sample_2002}). However, in many practical applications, the number of observations $n_i$ is small for some customers, making it challenging to make assumptions about the demand distribution in such cases, which prohibits the sampling process.

\begin{figure}
    \centering
    \includegraphics[scale=0.5]{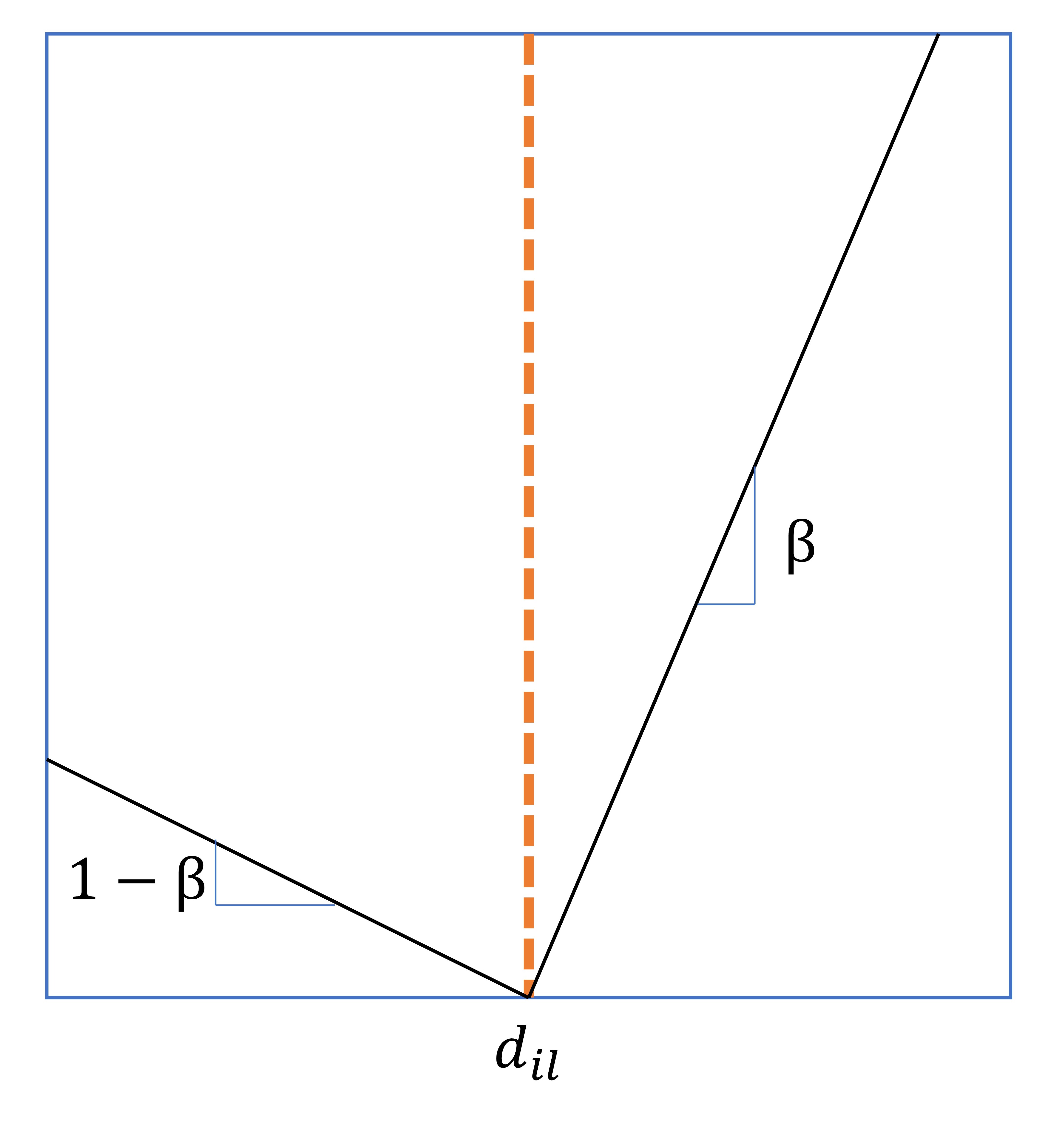}
    \caption{Tilted loss function}
    \label{fig:tilted_loss}
\end{figure}

Alternatively, the $\beta$ quantile of the demand can be estimated by minimizing the tilted loss (also called pinball loss) defined as 
\begin{equation} \label{eq:tilted_loss}
\frac{1}{n_i} \sum_{l=1}^{n_i} max(\beta(d-d_{il}),(\beta-1)(d-d_{il}))
\end{equation}
with respect to $d$ (see, for example, \cite{koenker_quantile_2001} for a more detailed discussion). Note that this quantile prediction requires no assumption on the distribution of demands. The tilted loss assigns asymmetric weight parameters to account for the under and over-prediction of quantiles. Figure \ref{fig:tilted_loss} shows how Equation \eqref{eq:tilted_loss} treats the larger and smaller value for $d_{il}$. However, as in the LSP case mentioned in the introduction, there might be new customers with no prior demand data, resulting in $n_i = 0$.

Although a small value of $n_i$ or $n_i=0$ prohibits individual demand predictions for each customer, there might be available feature information for all customers. In such cases, a contextual quantile prediction model can be employed to obtain demand predictions from customer features. This approach allows us to make predictions based on the contextual information shared by all customers, thereby overcoming the limitations of sparse demand data for specific individuals. 

Let $f_i \in \mathbb{R}^F$ represent the feature vector available for each customer $i$, where the feature information can include a combination of continuous values, binary values and categorical values. Then, we have the aggregated demand  history $\D = \{(f_i,d_{il}) : i \in \mathcal{C}, l \in \{1,\ldots,n_i\} \}$. In this case, it is possible to use a parametric function $q(f,\theta)$ as the predictor, where $f$ and $\theta$ are feature and parameter vectors, respectively. In this case, predicting the mean demand entails finding the optimal $\theta$ that minimizes the MSE loss function $\frac{1}{|\D|} \sum_{i} \sum_{l=1}^{n_i}(d_{il}-q(f_i,\theta))^2$.  The tilted loss function \eqref{eq:tilted_loss} can also be extended to the contextual case by defining
\begin{equation} \label{eq:ctx_tilted_loss}
\frac{1}{|\D|} \sum_{i} \sum_{l=1}^{n_i} max(\beta(q(f_i,\theta)-d_{il}),(\beta-1)(q(f_i,\theta)-d_{il}))
\end{equation}
as the loss function for our contextual quantile prediction model. The function $q$ can take various forms depending on the complexity desired in the prediction model. One straightforward option is to use an affine function, such as $q(f,\theta) = \theta_0 + \theta_1^Tf$. However, to capture nonlinear relationships, a more complex structure can be employed. For example, $q(f,\theta)$ can be defined as the output of a feed-forward neural network, with the parameters represented by $\theta$.  

Using the contextual quantile prediction model is also useful in order to predict uncertainty intervals used in RCVRPTW. For example, it is possible to predict $0.6$ and $0.95$ quantile values for $\overline{d}_i$ and $\overline{\overline{d}}_i$, respectively, and solve RCVRPTW with these predicted values. Note that occasionally, we can observe crossing quantiles, that is, $\overline{\overline{d}}_i < \overline{d}_i$ for one customer after the prediction. In this case, we simply let $\overline{\overline{d}}_i \leftarrow \max\{\overline{d}_i,\overline{\overline{d}}_i\}$ to avoid numerical problems in solving.  

\begin{figure}
    \centering
    \includegraphics[width=1\textwidth]{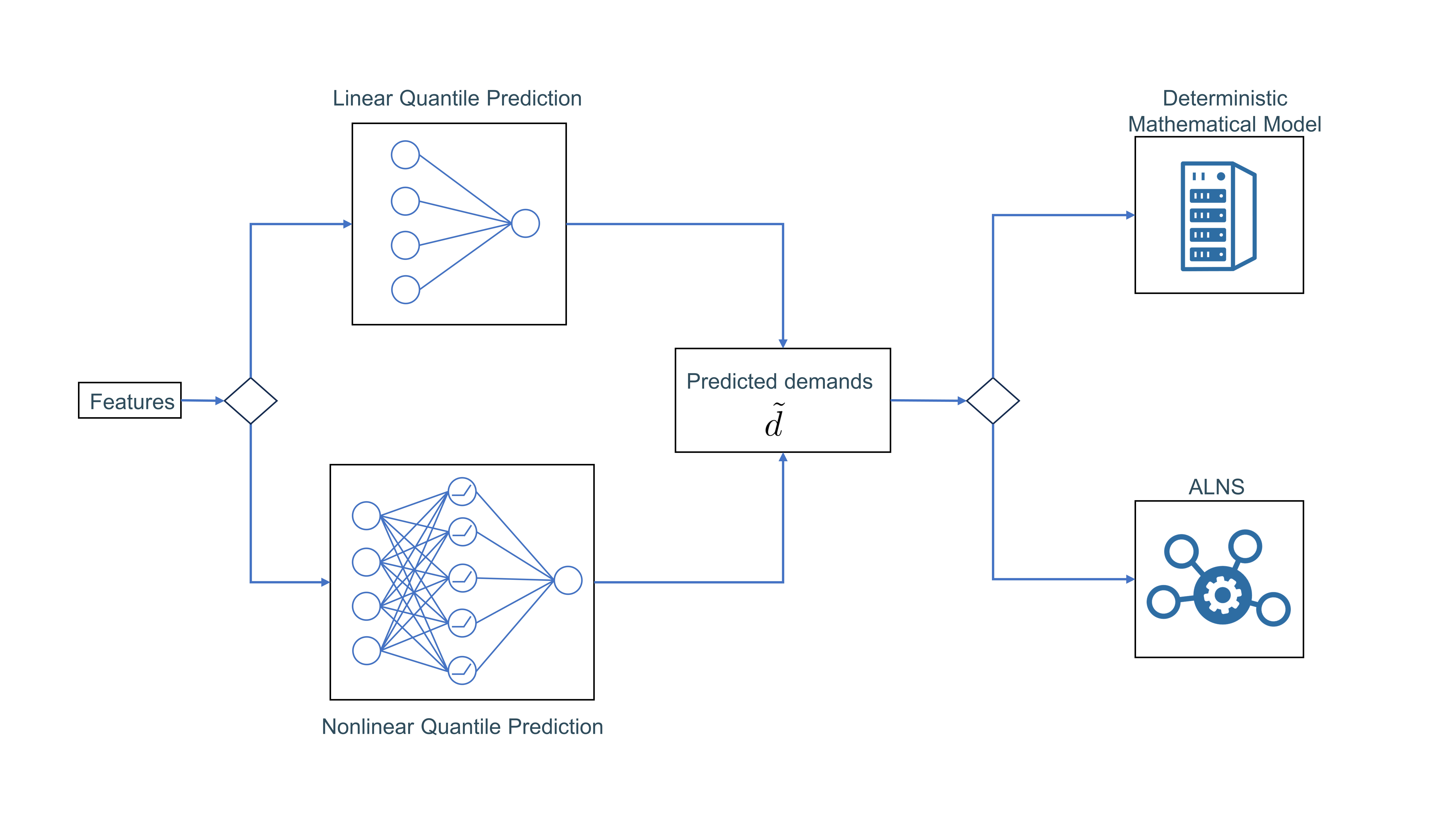}
    \caption{The complete process of decision-making for solving the deterministic model CVRPTW.}
    \label{fig:det_flow}
    \vspace{12pt}
    \centering
    \includegraphics[width=1\textwidth]{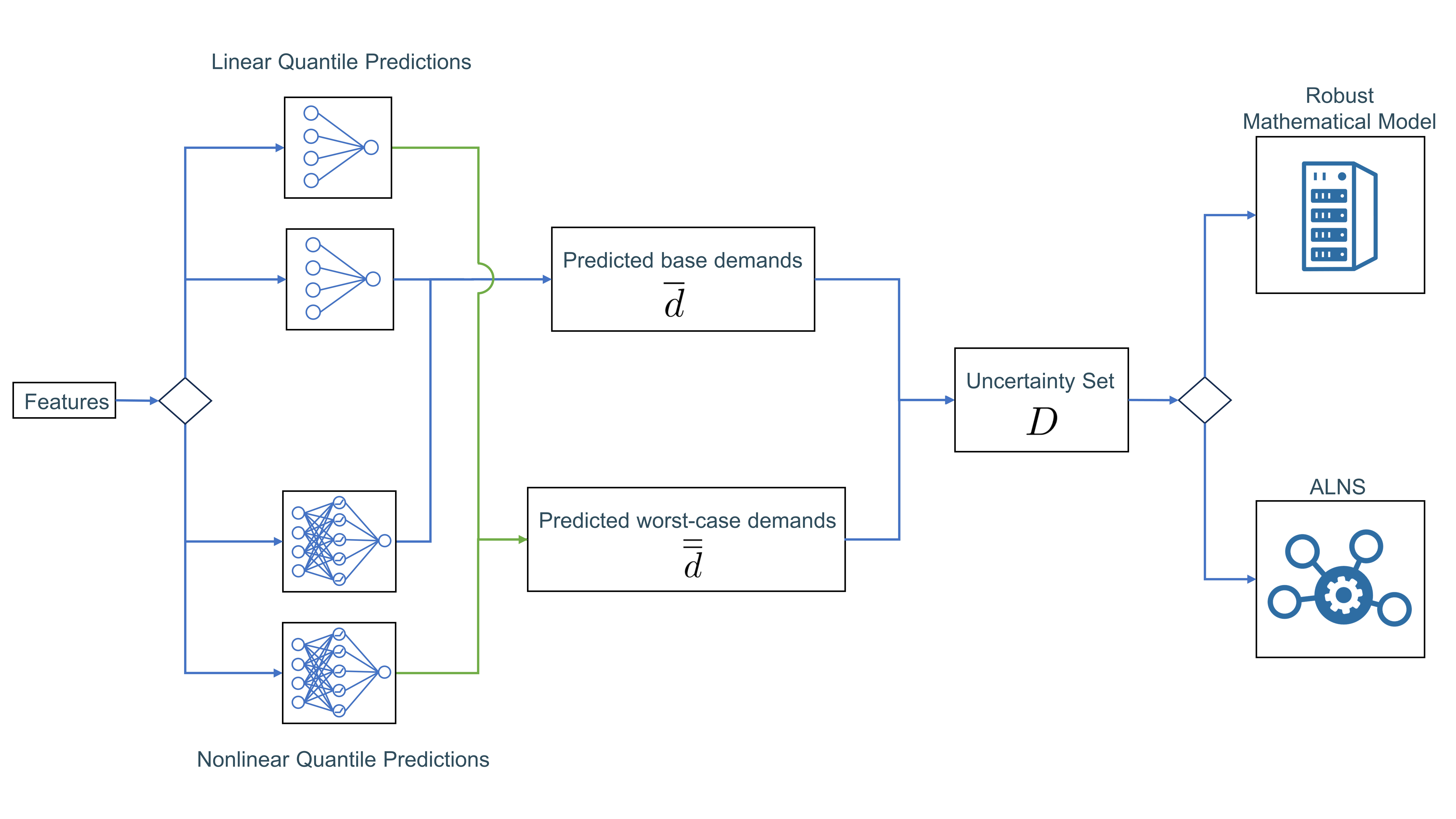}
    \caption{The complete process of decision-making for solving the robust model RCVRPTW.}
    \label{fig:rob_flow}
\end{figure}

Before presenting the results of the computational experiments, we provide a summary of the prediction and decision pipelines for CVRPTW and RCVRPTW in Figures \ref{fig:det_flow} and \ref{fig:rob_flow}, respectively. In both pipelines, the decision-maker first selects the method to predict demand quantiles (or means) based on observed features. For solving the deterministic problem, a single quantile value is predicted for each customer, while for the robust counterpart, two values representing the base and worst-case demand values are predicted. These predictions are used to define the uncertainty set $D$ in RCVRPTW. Following prediction, the decision-maker can opt for either an exact method or a heuristic approach to solve CVRPTW or RCVRPTW. 

\section{Computational Experiments}  \label{sec:experiments}

In this section, we first provide a description of the data and the experimental setup used in the computational experiments. Later, we present the results obtained for both the small (25-customer) and large (100-customer) instances.

For the experiments, we used a PC with 16 GB RAM and i9-11900H @ 2.50GHz processor. Python version 3.10.9 is used for the experiments and \cite{gurobi} version 10.0.1 is used to solve mixed integer programs. The details of implementation and the code to reproduce the experiments can be found at \url{https://github.com/irfanmahmutogullari/Leveraging-Contextual-Information-for-Robustness-in-Vehicle-Routing-Problems}.

\subsection{Experiment Setting} \label{sec:setting}

In our experiments, we use \cite{solomon_algorithms_1987} instances, a well-known benchmark in the literature.
The dataset includes 56 instances with 100 customers grouped into three categories C (c101-109,c201-208), R (r101-112,r201-211) and RC (rc101-108, rc201-208). To introduce customer features and randomness, we augment the dataset artificially. For each customer $i \in \mathcal{C}$, we create a 6-dimensional feature vector $f_i = (f_{i0},f_{i1},f_{i2},f_{i3},f_{i4},f_{i5})$ where $f_{i0}$ is the demand value of the customer in Solomon data. The other five features are created by assigning a random value between 0 and 1. The features $f_{i1}, f_{i2}, f_{i3}, f_{i4}$ and $f_{i5}$ could potentially reflect the economic and financial conditions of a customer, which can directly influence the customer's demand. The demand distribution $D_i$ of customer $i$ follows $D_i \sim \tilde{D}_i|\tilde{D}_i \geq 0$ where 

\begin{equation*}
    \tilde{D}_i \sim \begin{cases}
			f_{i0} + \mathcal{N}(-10,1) & \text{with probability } p_1\\
                f_{i0} + \mathcal{N}(-5,1)  & \text{with probability } p_2\\
                f_{i0} + \mathcal{N}(0,1)   & \text{with probability } p_3\\
                f_{i0} + \mathcal{N}(5,1)   & \text{with probability } p_4\\
                f_{i0} + \mathcal{N}(10,1)  & \text{with probability } p_5
		 \end{cases}
\end{equation*}
where $\mathcal{N}(\mu,\sigma)$ is the normal distribution with mean $\mu$ and standard deviation $\sigma$. The probabilities $p_1,\ldots,p_5$ are obtained by $p_k = \frac{e^{f_{ik}}}{\sum_{k' = 1}^{5}e^{f_{ik'}}}$ for $k \in \{1,\ldots,5\}$ and used to introduce multimodality and asymmetry in demand distribution.  

\begin{figure}
    \centering
    \includegraphics[scale=0.8]{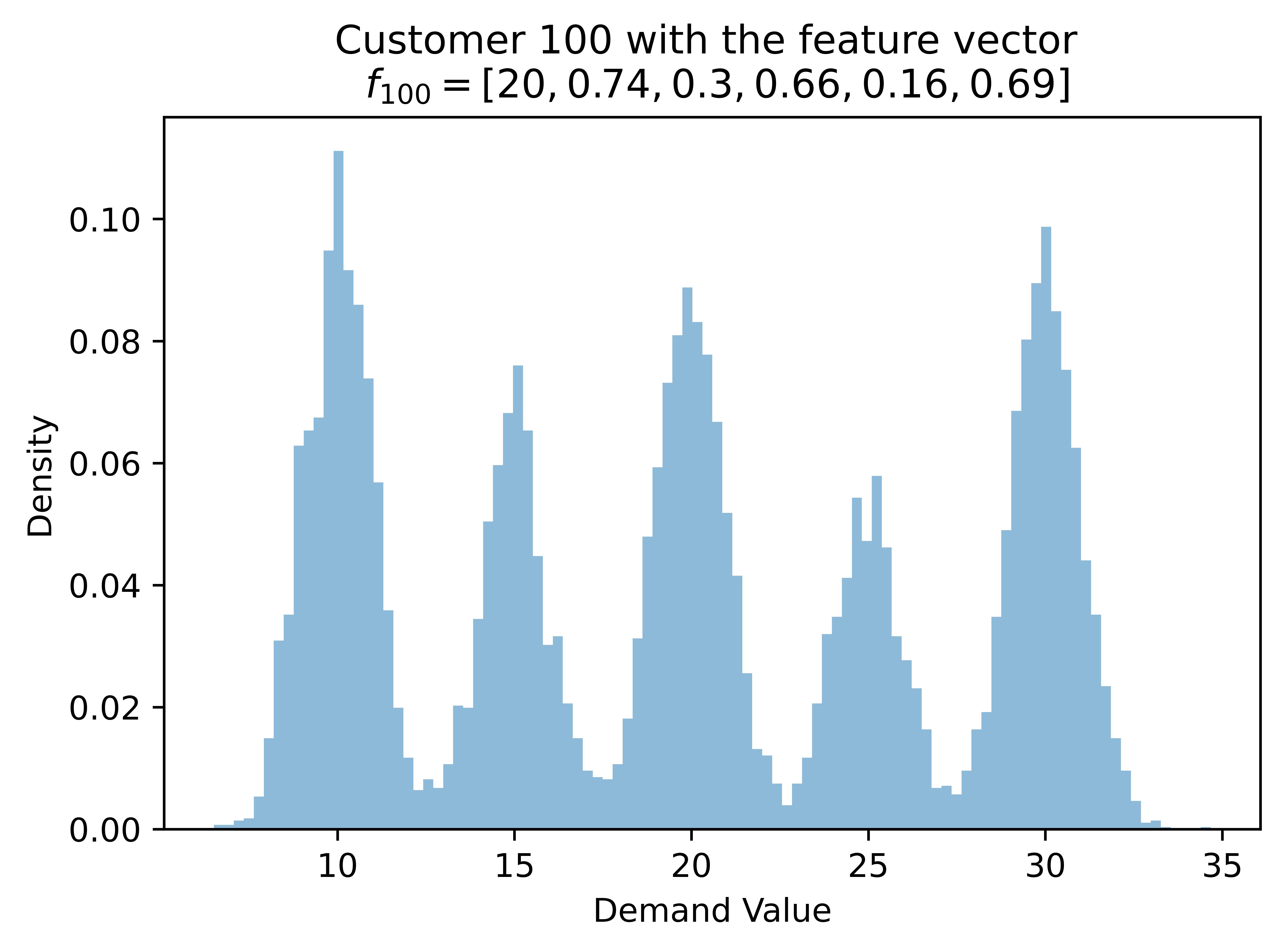}
    \caption{Histogram of demand for 10,000 sampled scenarios. (Instance = c101)}
    \label{fig:sample_histogram}
\end{figure}

Figure \ref{fig:sample_histogram} illustrates the empirical demand distribution of a customer in instance c101 for the given feature vector as an example. We also halved the vehicle capacities in the original data set to observe the effect of demand uncertainty more clearly. 

In our experiments, we investigated various settings concerning the demand history. We conducted tests for scenarios where we have demand history for all 100 customers (\textit{all}), half of the customers (\textit{half}), and a quarter of the customers (\textit{quar}). For customers with available demand history, the number of observations ($n_i$) in the data history can take one of three values: 1, 10, and 30. As a result, we have a total of nine different settings, each representing a different level of data.

After obtaining the demand data, a prediction model is employed to forecast the demand values to be used in the decision-making process. When not using contextual information, individual (I) mean or quantile predictions can be made for customer $i \in \mathcal{C}$ if $n_i > 1$. On the other hand, in the contextual setting, both linear (L) and nonlinear (N) prediction models can be used for predicting the mean and quantiles. In our experiments, we use a straightforward feedforward neural network with a single hidden layer of 10 neurons with the ReLU activation function. The feature values are normalized during the training. In order to create and train prediction models, we used the Keras library. The maximum number of iterations in a training is set to 1,000 and the training is stopped if no improvement is made in the last 50 iterations. 

\begin{figure}
    \centering
    \includegraphics[scale=0.9]{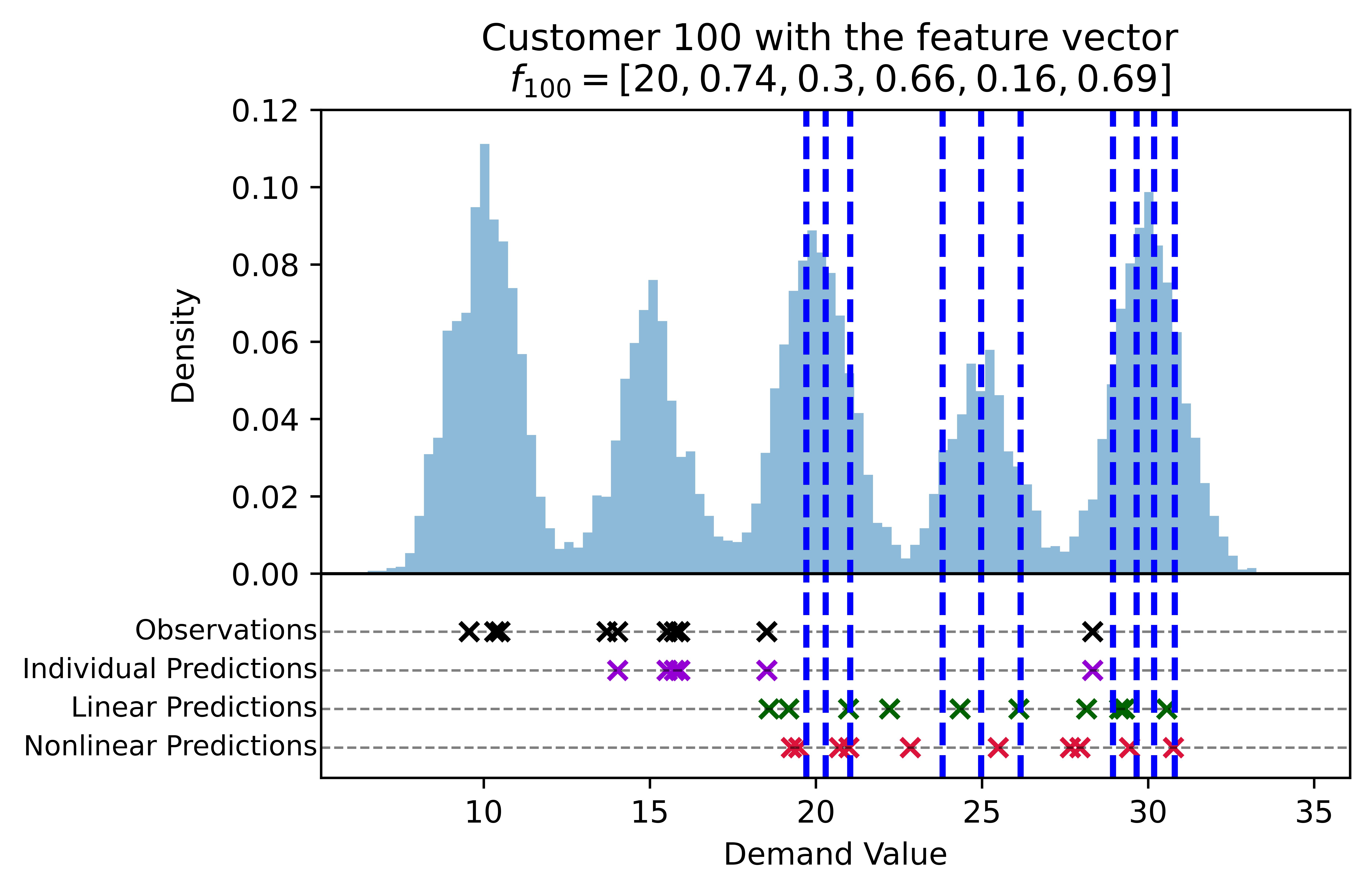}
    \caption{The blue vertical lines represent the empirical quantile predictions obtained from sampling 10,000 scenarios. The black ``x" represents the observed demand values. The violet, green and red ``x" represent individual, linear and nonlinear quantile predictions (Instance = c101).}
    \label{fig:sample_histogram_with_quantiles}
\end{figure}

Figure \ref{fig:sample_histogram_with_quantiles} illustrates the quantile predictions for the example customer mentioned in Figure \ref{fig:sample_histogram}. The blue vertical lines correspond to the empirical 0.50, 0.55,\ldots,0.95 quantile predictions obtained from sampling 10,000 scenarios. When there are 10 observed demand values (indicated by black ``x"), the decision maker can make non-contextual quantile predictions (indicated by violet ``x") using the realized demands. However, these predictions are highly influenced by the observed data and may not provide accurate quantile estimates. On the other hand, the contextual quantile linear and nonlinear predictions (indicated by a green and red ``x", respectively) offer improved predictions that better align with the actual quantile values.
 
Using the predicted demand values, the decision maker can solve a deterministic CVRPTW (D) or its robust counterpart RCVRPTW (R) to make routing decisions. For small instances, it is feasible to solve the optimization models presented in Section \ref{sec:models} optimally using a commercial solver. However, for larger instances, finding the exact solution may not be possible within reasonable time limits. In such cases, the heuristic presented in Section \ref{sec:alns} can be employed to obtain near-optimal solutions efficiently.

For the deterministic models, we adopt the naming convention D-$\delta$-$\epsilon$, where D refers to \textit{deterministic}, $\delta \in \{\text{I (individual), L (linear prediction model), N (nonlinear prediction model)}\}$ indicates the prediction method and $\epsilon \in \{\text{M (indicating mean)},50,55,\ldots,90,95\}$ represents the mean or quantile prediction utilized in the deterministic model. For instance, D-I-M indicates the solution of the deterministic model, where the mean demand values are individually used for each customer in the non-contextual setting. On the other hand, D-N-60 indicates that the 60th quantile predictions are obtained using the nonlinear prediction model in the deterministic problem.

Similarly, for the robust setting, we employ the naming convention R-$\delta$-$\overline{\epsilon}$-$\hat{\epsilon}$-$\Gamma\gamma$, where R denotes \textit{robust} and $\delta \in \{\text{I, L, N}\}$ indicates the prediction method. $\overline{\epsilon}$ and $\hat{\epsilon}$ indicate the levels of $\overline{d}$ and $\hat{d}$, respectively, in the robust setting. Finally, $\gamma$ denotes the budget parameter in the definition of the uncertainty set \eqref{eq:polyhedral_set}. For example, R-L-M-95-$\Gamma$1 indicates that the mean and 95th percentile predictions obtained by the linear prediction model are used for the robust problem with $\Gamma = 1$.  

To evaluate the quality of the obtained routes, a simulation is performed by sampling 10,000 demand values from the true distribution for each customer. The total cost for these routes is calculated, which includes the initial cost and the recovery cost. The initial cost of a route is the sum of the used arc costs as given in \eqref{eq:det-obj} and \eqref{eq:rob-obj}. However, due to the uncertainty of demand values, some routes may become infeasible after demand realization. In such cases, the recovery cost is the total cost induced by the detours. We assume that customer $i$ is served at the earliest possible time, which is the minimum of $T^{min}_i$ and the reach time to customer $i$. If the service beginning time of the customer exceeds $T^{max}_i$, a time window violation occurs, indicating that the customer cannot be served within the specified time window.  Note that the initial routes satisfy time window constraints for the predicted demand vector; however, detours caused by uncertain demand can lead to potential violations.  

\subsection{Results for Small Instances}

We first present the results of the experiments for the C instances with 25 customers. The objective of these experiments is to observe the performance of the proposed setting when an exact solver is used to solve the MIP formulation of the problem. By analyzing the results for these smaller instances, we can gain insights and make informed decisions when scaling up to larger problem instances.

For each Solomon instance, we consider the first 25 customers in the problem. In total, there are 17 instances and for each instance, we considered 9 different settings for the available data. In this settings, \textit{all} indicates the demand history exists for all 100 customer. On the other hand, for \textit{half} and \textit{quar}, the demand history is available for the last 50 and 25 customers in the original instance. 

For each instance and setting, the following models are solved \begin{itemize}
\item D-$\delta$-$\epsilon$ for all $\delta \in$ \{I,L,N\}, $\epsilon \in \{M,50,55,\ldots,90,95\}$

\item R-$\delta$-$\overline{\epsilon}$-$\hat{\epsilon}$-$\Gamma\gamma$ for all $\delta \in$ \{I,L,N\}, $\overline{\epsilon} \in \{M,50,55\}$, $\hat{\epsilon} \in \{90,95\}$, $\gamma \in \{1,2\}$
\end{itemize} 
for 5 randomly generated instances.  (Note that, for \textit{half} and \textit{quar} the customers considered in the problem do not have a demand history so the models D-I-$\epsilon$ and R-I-$\overline{\epsilon}$-$\hat{\epsilon}$-$\Gamma\gamma$ are not available.) In total, there are 39,100 instances. Among these, we applied a time limit of 60 seconds for each instance, and 20,322 instances were solved optimally within this time constraint. For the remaining instances, we kept the best solution obtained within the time limit.

\begin{figure}[htbp]
    \centering
    
    \begin{subfigure}{0.332\textwidth}
        \centering
        \includegraphics[width=\linewidth]{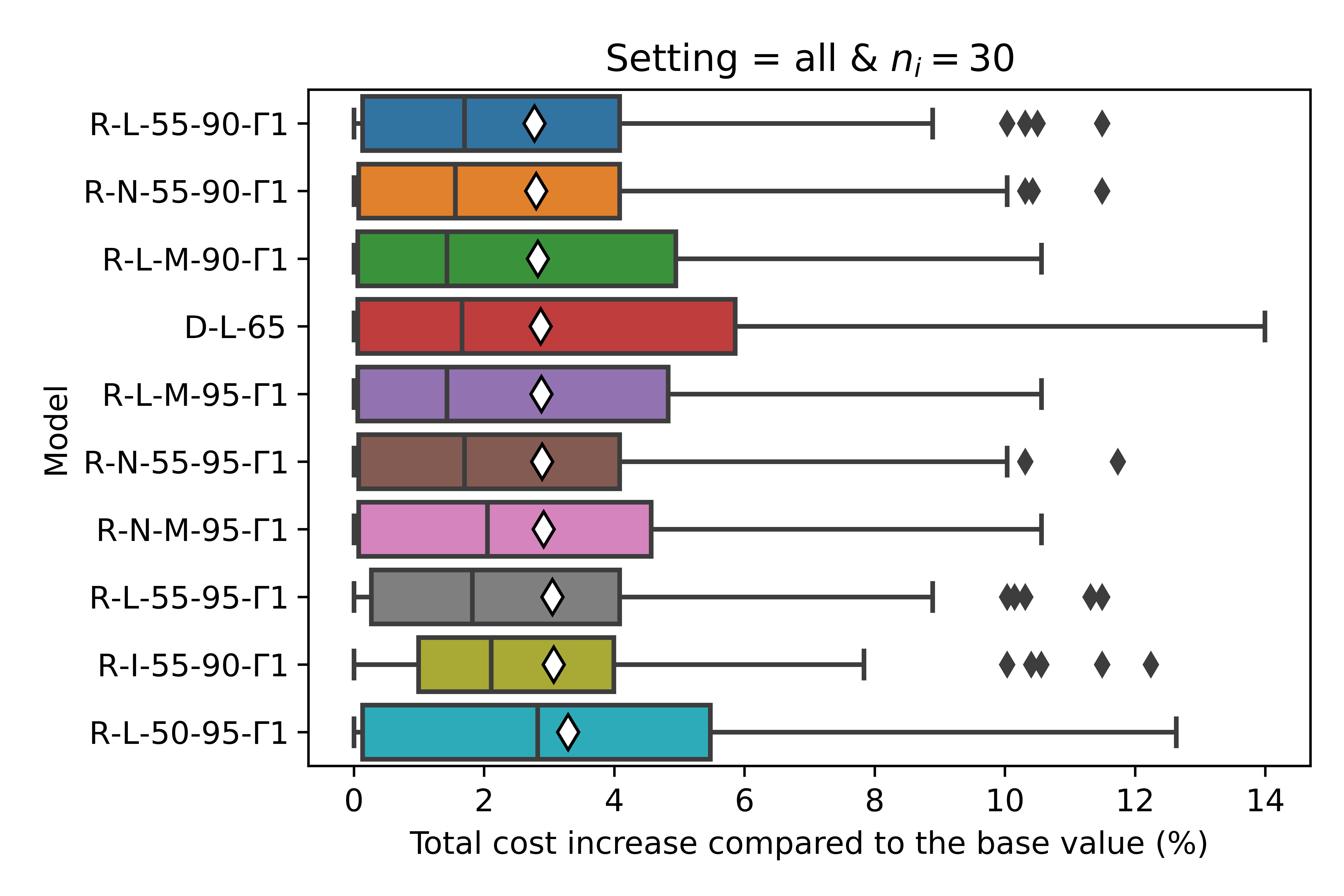}
    \end{subfigure} \hspace{-0.2cm} \begin{subfigure}{0.332\textwidth}
        \centering
        \includegraphics[width=\linewidth]{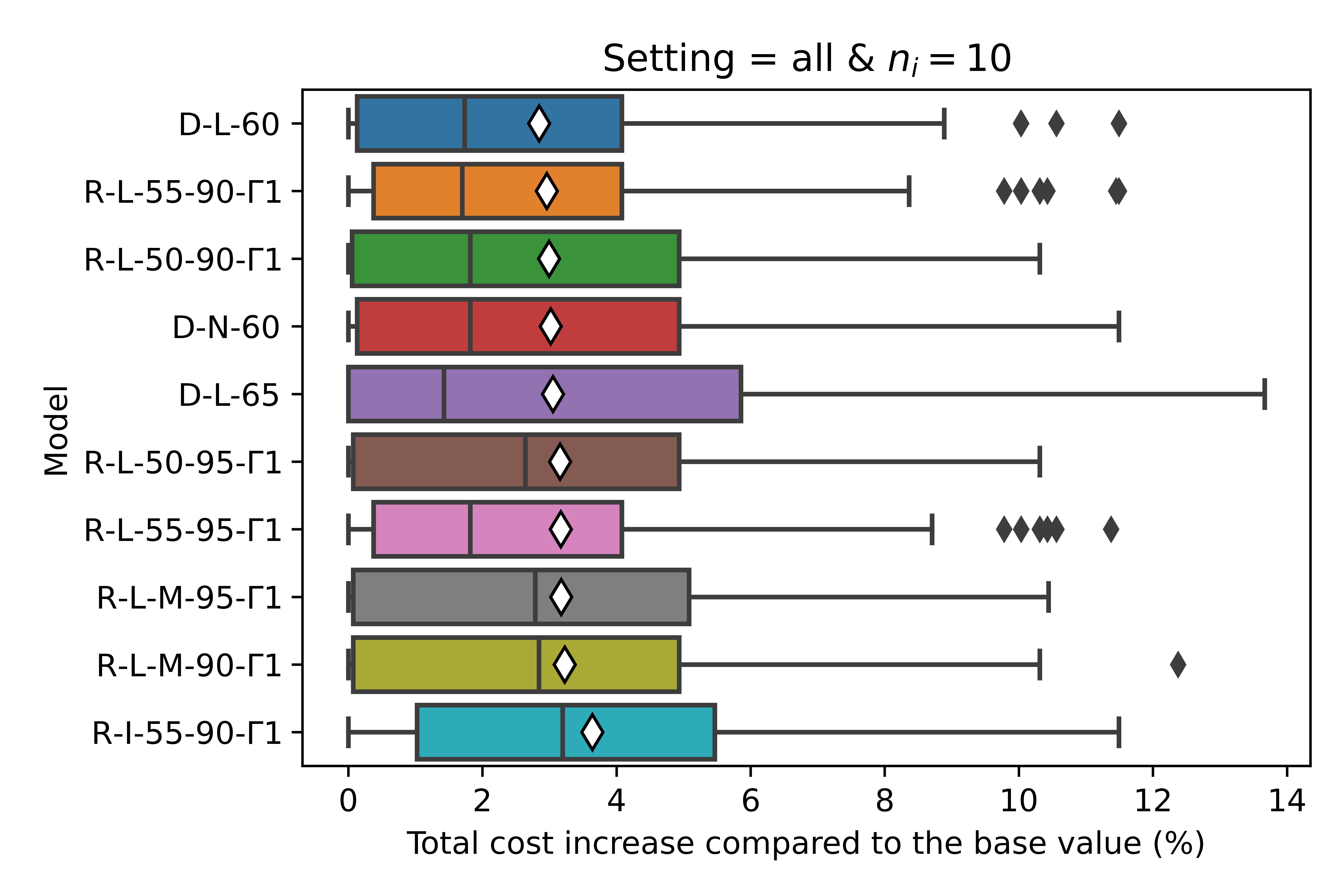}
   \end{subfigure} \hspace{-0.2cm}  \begin{subfigure}{0.332\textwidth}
        \centering
        \includegraphics[width=\linewidth]{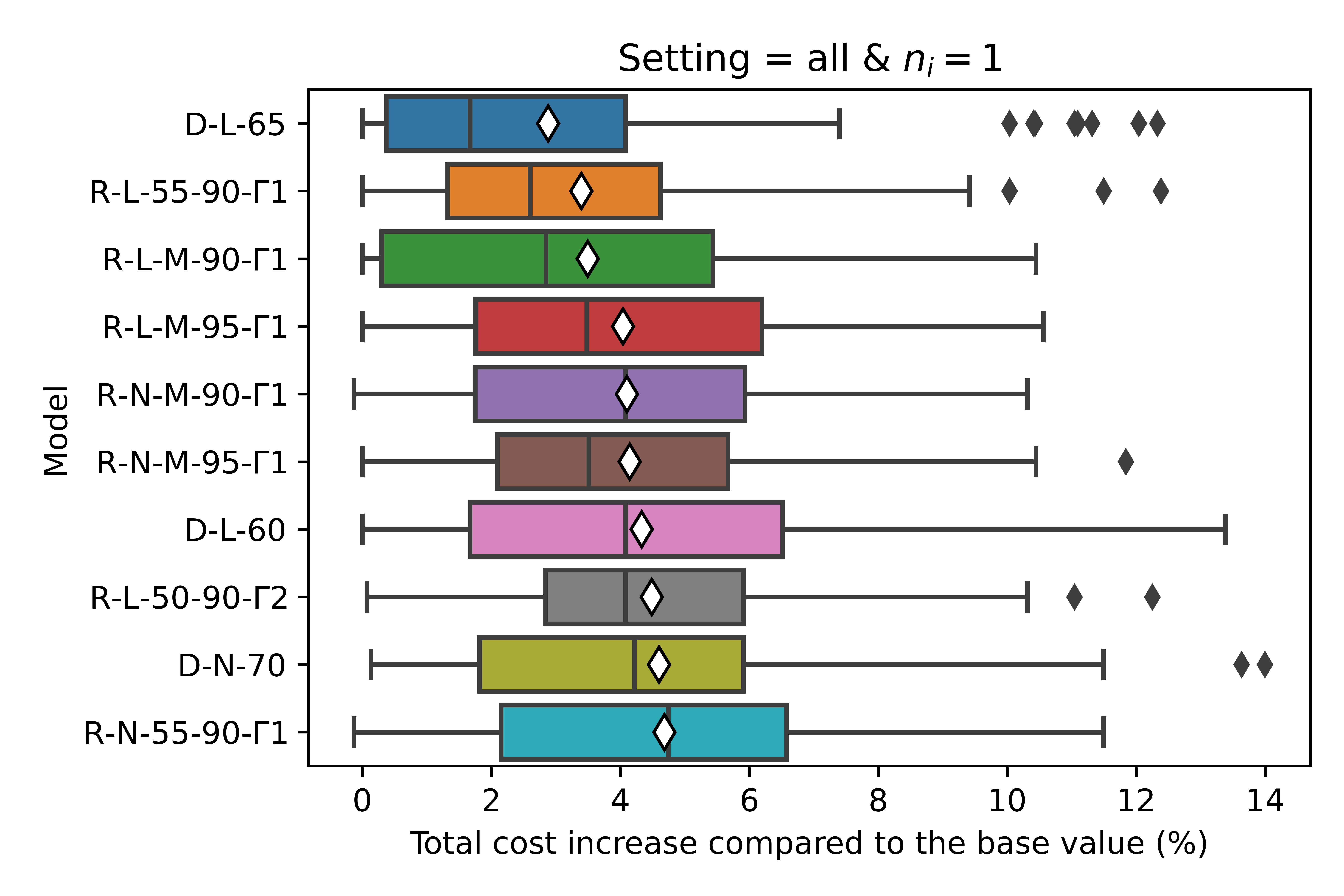}
   \end{subfigure}

    \vspace{0.5cm} % Add some vertical space between rows
    
    \begin{subfigure}{0.332\textwidth}
        \centering
        \includegraphics[width=\linewidth]{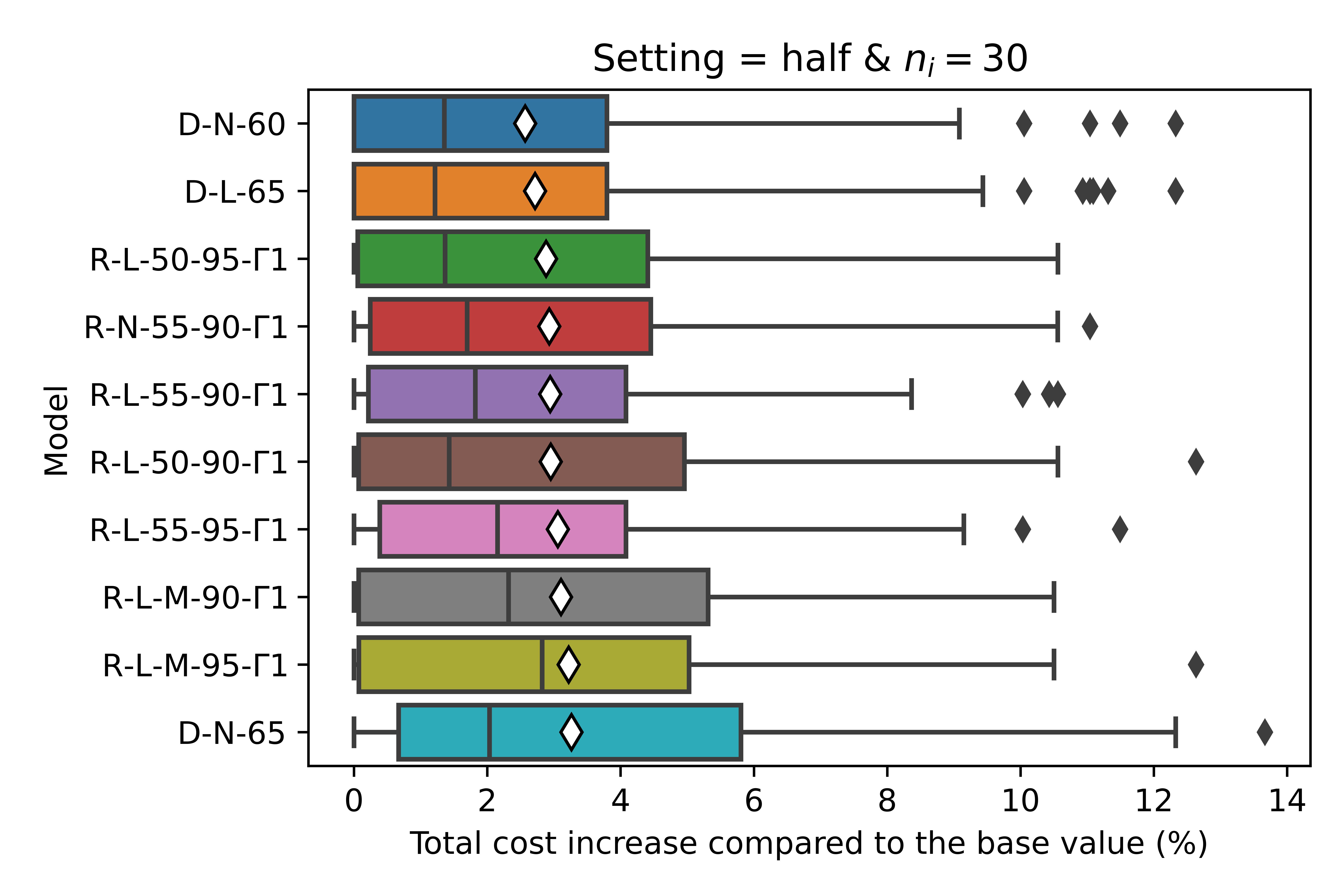}
           
    \end{subfigure}  \hspace{-0.2cm} \begin{subfigure}{0.332\textwidth}
        \centering
        \includegraphics[width=\linewidth]{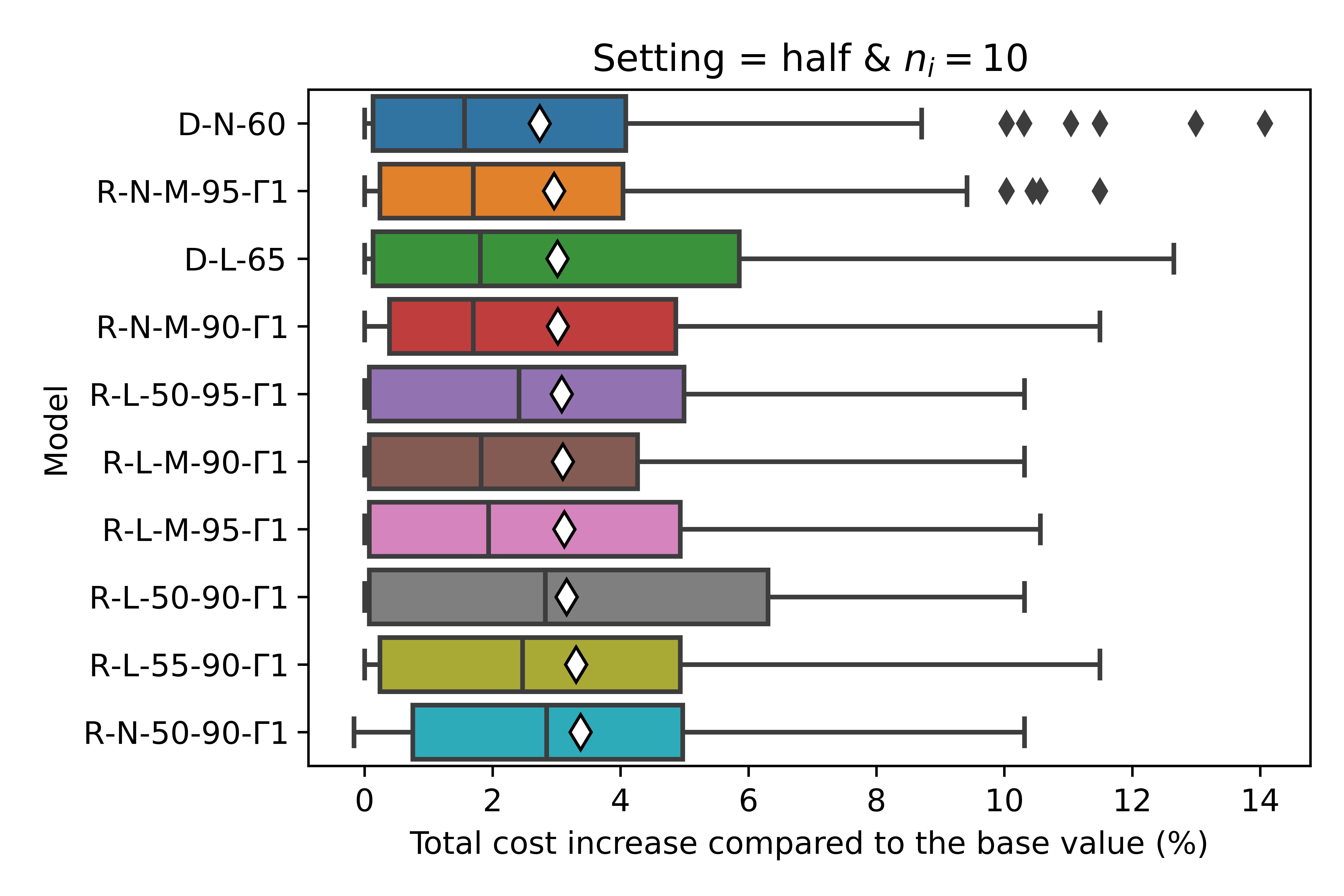}
         
    \end{subfigure} \hspace{-0.2cm}  \begin{subfigure}{0.332\textwidth}
        \centering
        \includegraphics[width=\linewidth]{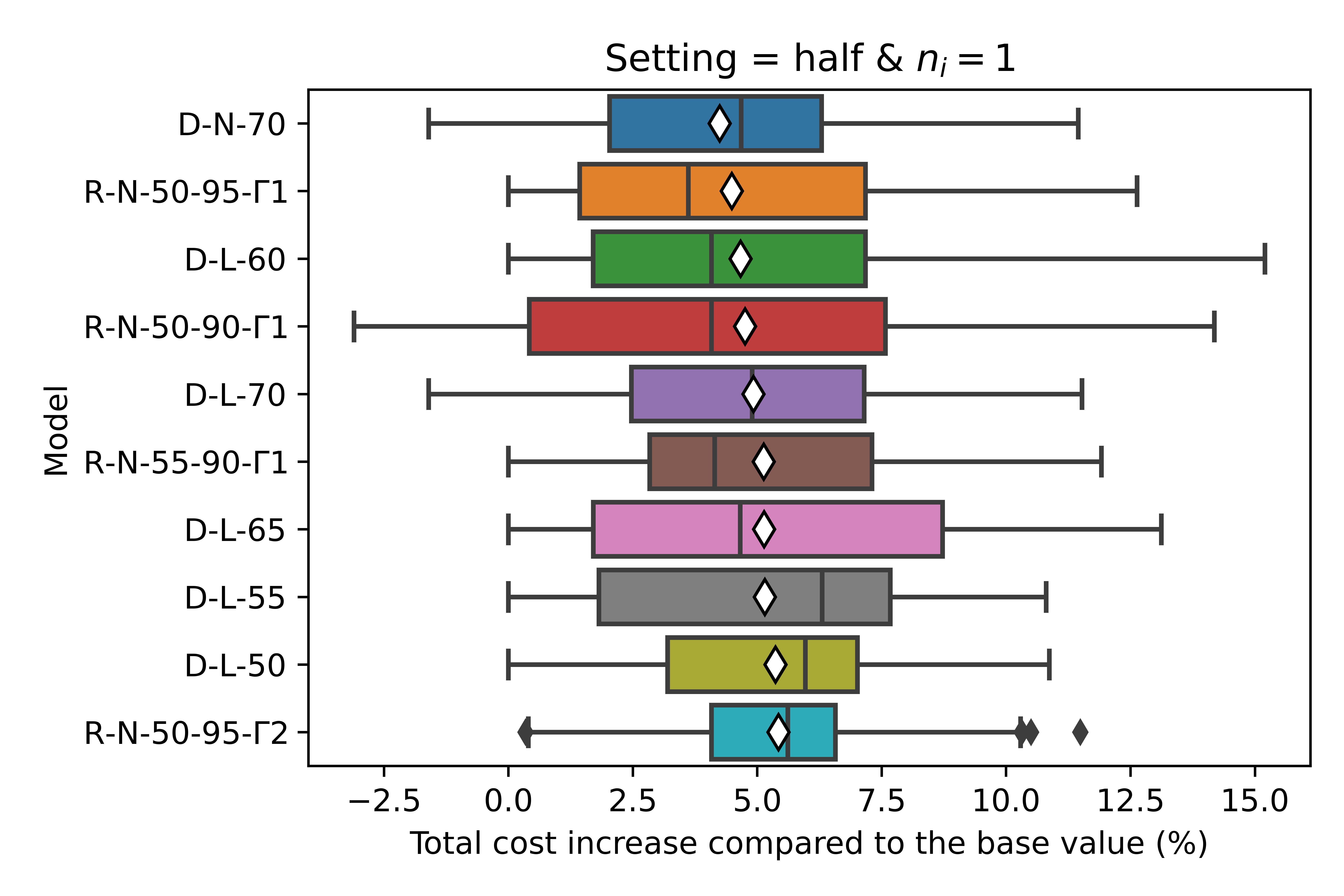}
    \end{subfigure}
    
    \vspace{0.5cm} % Add some vertical space between rows
    
    \begin{subfigure}{0.332\textwidth}
        \centering
        \includegraphics[width=\linewidth]{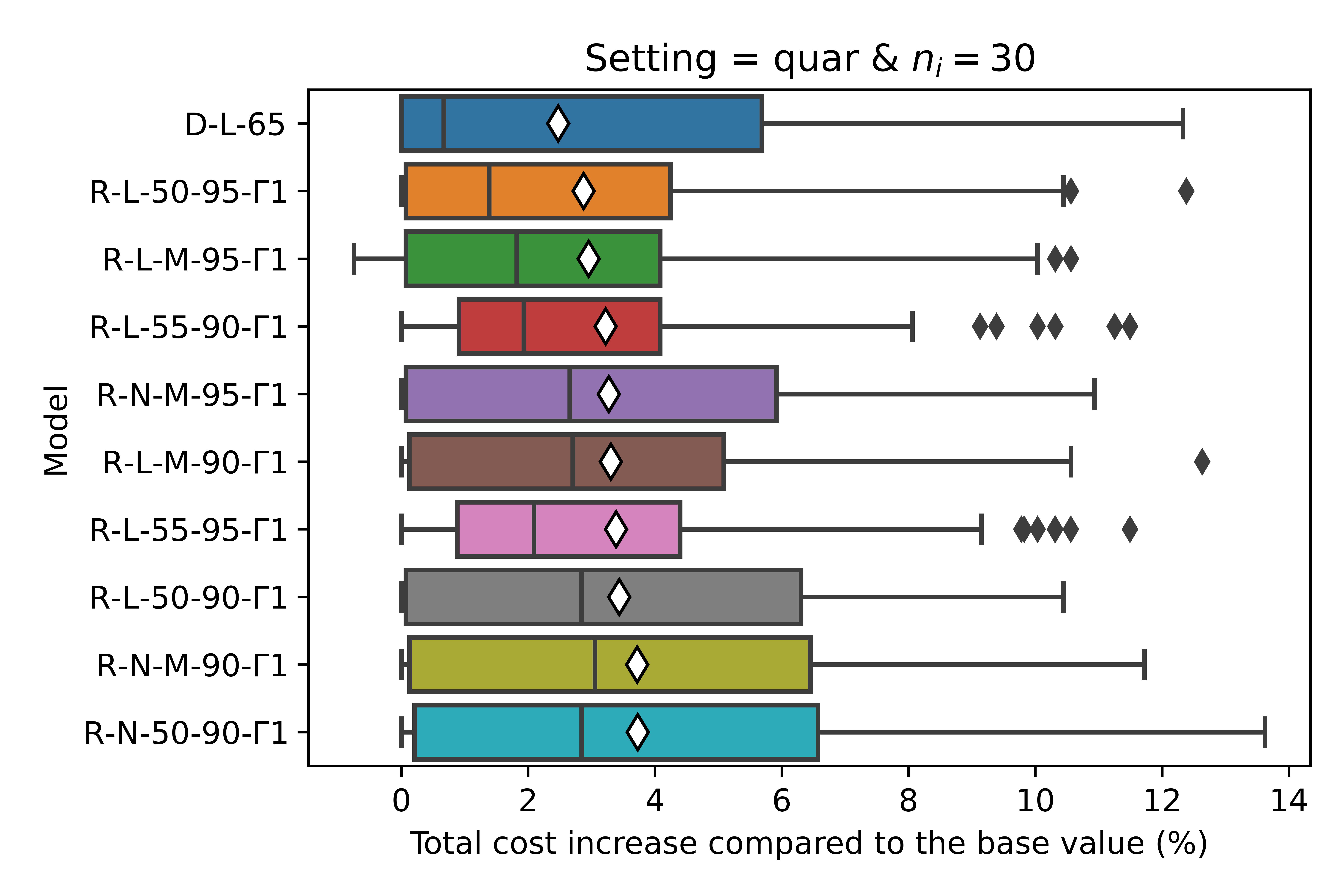}
          
    \end{subfigure} \hspace{-0.2cm} \begin{subfigure}{0.332\textwidth}
        \centering
        \includegraphics[width=\linewidth]{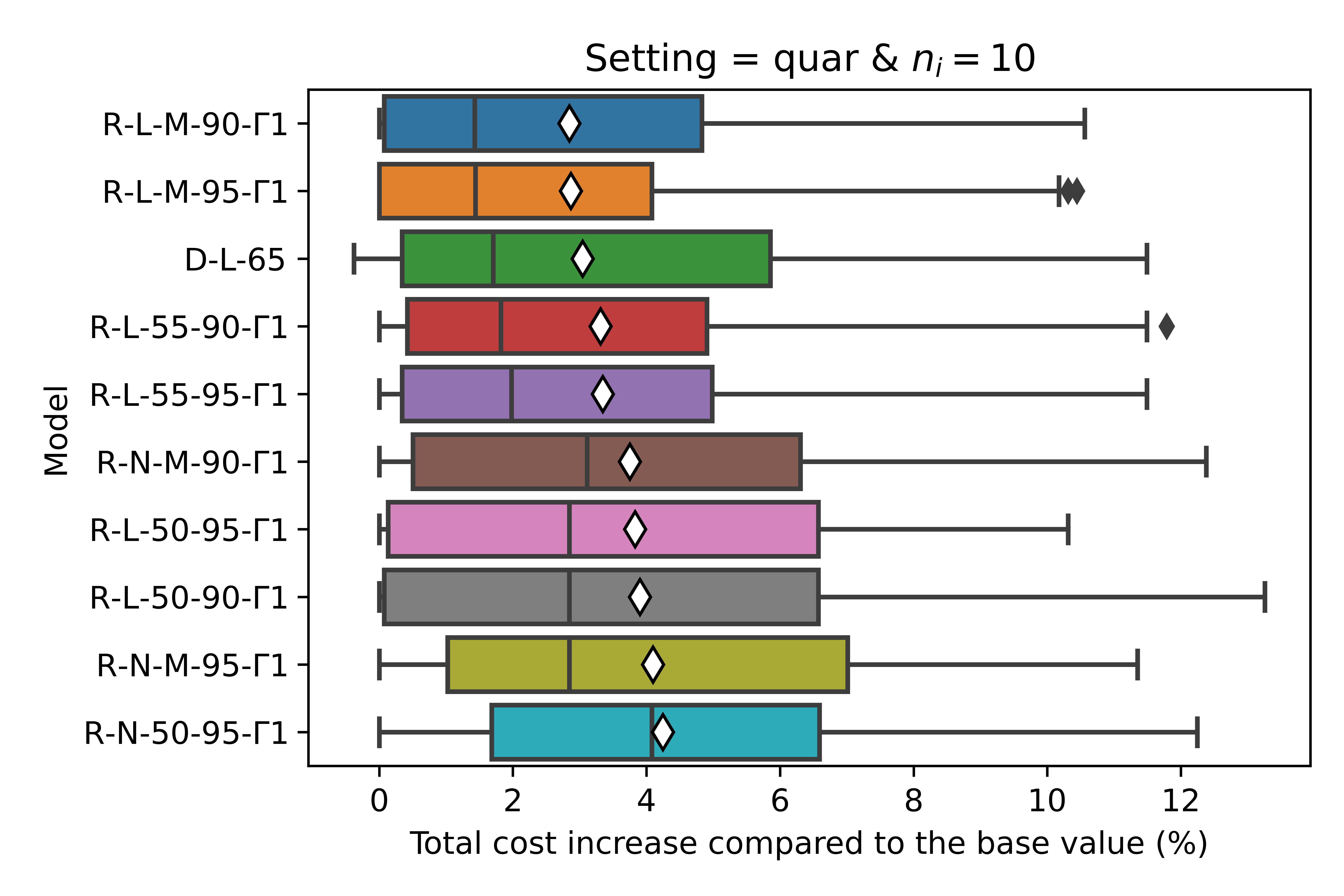}
         
    \end{subfigure}  \hspace{-0.2cm} \begin{subfigure}{0.332\textwidth}
        \centering
        \includegraphics[width=\linewidth]{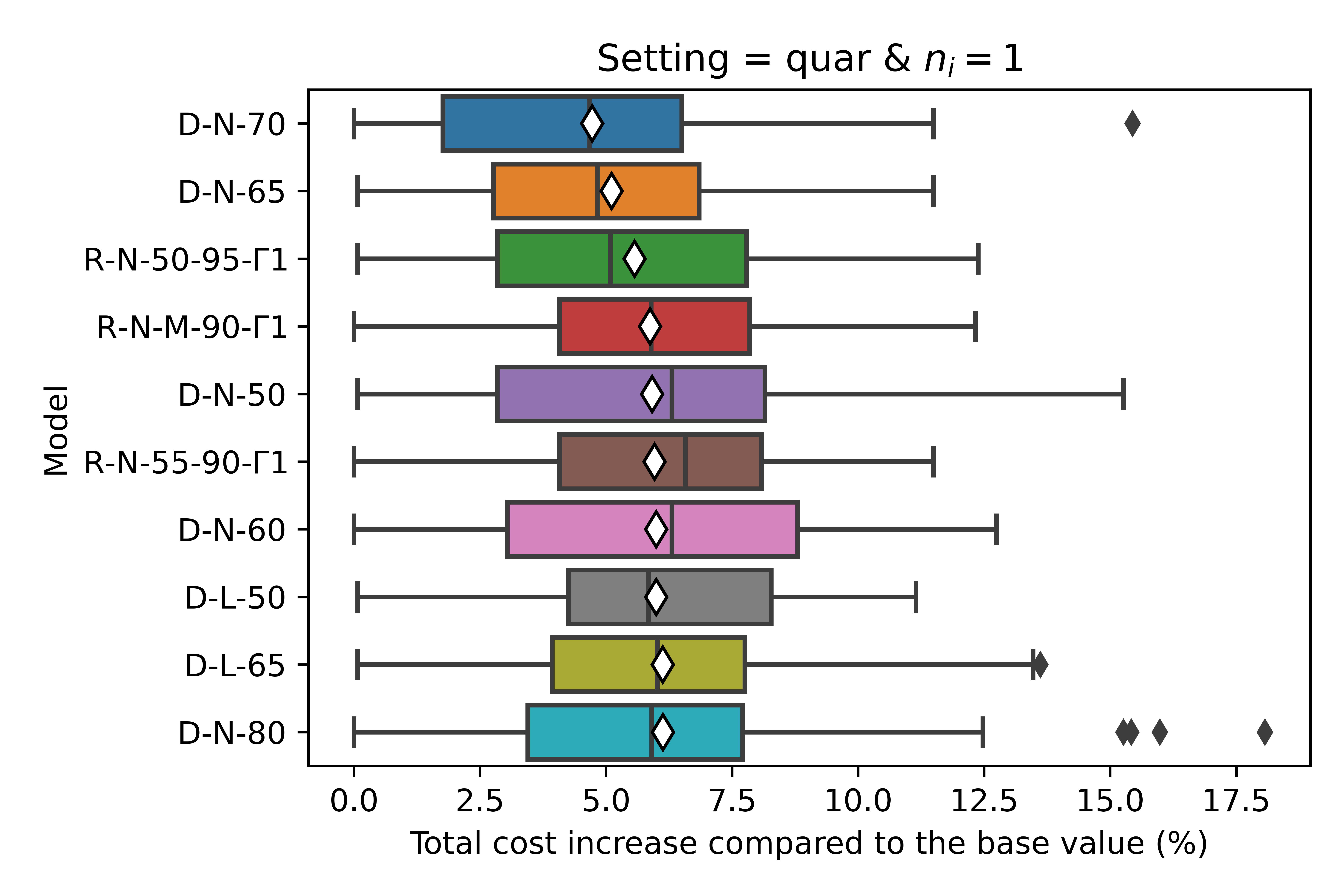}
          
    \end{subfigure}
    \caption{Performances of models for different settings instances with 25 customers. White and black diamonds represent the average values and outliers, respectively.}
    \label{fig:results_small_instances}
\end{figure}

Figure \ref{fig:results_small_instances} shows the normalized total costs obtained from the simulations for different models and settings. For each setting, only the 10 best-performing models are presented. The costs are normalized by calculating their percentage deviation from the base value for each instance. The base value is determined as the best performance achieved by any model in the \textit{all} setting with $n_i$ = 30. It is worth noting that this value can be negative, as a model from another setting may occasionally outperform the best model in the \textit{all} setting with $n_i$ = 30.
 
As depicted in Figure \ref{fig:results_small_instances}, the most data-rich scenario, specifically setting \textit{all} with $n_i$ = 30, shows that the robust model performs best on average. The model R-L-55-90-$\Gamma$1 exhibits an average cost increase of 2.77\% (with a range of 0\% to 11.49\%). Additionally, 9 out of the 10 best-performing models are robust models, suggesting that in a data-rich environment, predicting multiple quantiles and employing a robust model for decision-making is a promising approach. Nevertheless, the deterministic model D-L-65 also shows a competitive performance with a mean cost increase of 2.86\% (with a range of 0\% to 13.99\%), indicating that even the deterministic model can provide some level of robustness by using a value 65th quantile greater than the mean or median prediction.

Another observation from Figure \ref{fig:results_small_instances} is the impact of the amount of available information on the decision quality. As we move from \textit{all} to \textit{quar}, and from $n_i = 30$ to $n_i = 1$, we observe an increase in the average total cost. This outcome is intuitive, as having more data allows for better prediction, leading to improved performance in the decision-making process.

As mentioned in the previous section, our aim is to achieve robustness in our decisions. Interestingly, we observe that the solutions of the deterministic models can also lead to robust decisions, especially in data-scarce settings. For instance, when the setting is \textit{quar} with $n_i=1$, the solution of D-N-70 gives the most robust result with an average cost increase of 4.72\% (ranging from 0\% to 15.44\%). Moreover, 7 out of the 10 best-performing models in terms of total cost are deterministic models. This suggests that a more complex prediction strategy, such as predicting multiple quantiles accurately, often fails due to data scarcity. It is noteworthy that this result aligns with \cite{elmachtoub_estimate-then-optimize_2023}, which demonstrates that in cases with sufficient available data, complex prediction models yield better performance in the subsequent decision-making process. 

Another advantage of deterministic models is that they generally require simpler models and shorter solution times. Table \ref{tab:model_runtimes} presents the average and standard deviation of the running times for both the deterministic and robust models. On average, the running time for the deterministic model is smaller than that of the robust models. Additionally, when the level of robustness is set high ($\Gamma = 2$), the robust model requires even more time compared to the low level of robustness $\Gamma = 1$. Additionally, 56.25\% of the deterministic instances were solved optimally within the time limit. For the cases where $\Gamma$ equals 1 or 2, this fraction is slightly lower at 51.89\% and 47.17\%, respectively.

\begin{table}[h]
  \centering
  \caption{The mean and standard deviation of running times (in seconds) for the deterministic and robust model and the fraction of instances solved optimally within the time limit.}
    \begin{tabular}{cccc}
          & \multicolumn{2}{c}{Solution Time } &  Fraction of Instances \\
\cmidrule{2-3}    Model & Mean  & St. Deviation  & Solved Optimally  \\
    \midrule
    Deterministic & 33.94 & 25.35  &56.25 \%\\
    Robust $\Gamma$ = 1 & 35.18 & 26.01  &51.89 \%\\
    Robust $\Gamma$ = 2 & 38.53 & 24.35  &47.17 \%\\
    \end{tabular}%
  \label{tab:model_runtimes}%
\end{table}%

The final observation from Figure \ref{fig:results_small_instances} is that contextual prediction strongly outperforms non-contextual individual predictions in terms of robustness. Only in the most data-rich setting, \textit{all} with $n_i=30$, R-I-55-90-$\Gamma$1 appears among the 10 best-performing models. In all other cases, the contextual prediction models demonstrate superior robustness compared to individual predictions. This highlights the importance of utilizing contextual information when making demand predictions, especially in scenarios with limited data for individual customers.

Another analysis can be conducted regarding the violation of customer time windows. As anticipated, the smallest number of time window violations can be achieved with very conservative models such as D-$\delta$-95 or R-$\delta$-$55$-$95$-$\Gamma2$. The reason behind this is that these models are overly pessimistic, making their route decisions based on excessively predicted demand values. While this approach ensures that the resulting solution is feasible for most demand realizations, it may not be cost-effective. 

Therefore, in Table \ref{tab:small_instances_TW}, we present the relative increase in the fraction of customers with time window violations for the models presented in Figure \ref{fig:results_small_instances} with respect to the model D-$\delta$-95, $\delta \in $ \{I,L,N\} which gives the lowest time window violation. 
 

\begin{table}[htbp]
  \centering
  \caption{Average percentage increase in the customers with time window violations compared to the most conservative model D-$\delta$-95 for each setting and number of observations ($n_i$) for small instances.}
    \begin{tabular}{c|cc|cc|cc}
    \multicolumn{1}{c}{} & \multicolumn{6}{c}{Setting} \\
\cmidrule{2-7}    $n_i$  & \multicolumn{2}{c}{\textit{all}} & \multicolumn{2}{c}{\textit{half}} & \multicolumn{2}{c}{\textit{quar}} \\
    \midrule
    \midrule
    30    & Model  & Mean (Std) & Model  & Mean (Std) & Model  & Mean (Std) \\
    \midrule
          & D-L-65 & 9.76 (9.33) & D-N-65 & 10.11 (9.48) & R-N-M-90-$\Gamma$1 & 9.39 (9.28) \\
          & R-L-50-95-$\Gamma$1 & 10.36 (9.41) & R-L-55-90-$\Gamma$1 & 10.57 (8.92) & R-L-50-90-$\Gamma$1 & 9.87 (8.79) \\
          & R-L-55-90-$\Gamma$1 & 10.53 (9.94) & R-L-M-90-$\Gamma$1 & 10.59 (8.63) & D-L-65 & 9.9 (9.34) \\
          & R-N-55-90-$\Gamma$1 & 10.8 (8.79) & R-N-55-90-$\Gamma$1 & 10.68 (8.63) & R-N-M-95-$\Gamma$1 & 10.22 (9.3) \\
          & R-N-55-95-$\Gamma$1 & 10.81 (8.77) & R-L-M-95-$\Gamma$1 & 10.88 (9.67) & R-L-M-95-$\Gamma$1 & 10.57 (8.54) \\
          & R-L-M-95-$\Gamma$1 & 11.28 (9.76) & R-L-50-95-$\Gamma$1 & 10.9 (9.53) & R-N-50-90-$\Gamma$1 & 10.86 (9.45) \\
          & R-N-M-95-$\Gamma$1 & 11.28 (9.76) & R-L-55-95-$\Gamma$1 & 11.05 (9.37) & R-L-50-95-$\Gamma$1 & 10.87 (9.64) \\
          & R-I-55-90-$\Gamma$1 & 11.32 (9.3) & D-N-60 & 11.4 (8.53) & R-L-55-95-$\Gamma$1 & 10.9 (8.97) \\
          & R-L-55-95-$\Gamma$1 & 11.56 (9.61) & D-L-65 & 11.65 (9.54) & R-L-M-90-$\Gamma$1 & 11.06 (8.61) \\
          & R-L-M-90-$\Gamma$1 & 11.78 (9.46) & R-L-50-90-$\Gamma$1 & 11.84 (10.26) & R-L-55-90-$\Gamma$1 & 11.64 (8.39) \\
    \midrule
    10    & Model  & Mean (Std) & Model  & Mean (Std) & Model  & Mean (Std) \\
    \midrule
          & R-I-55-90-$\Gamma$1 & 11.43 (9.56) & D-L-65 & 8.94 (9.48) & R-L-50-90-$\Gamma$1 & 8.73 (9.26) \\
          & R-L-M-95-$\Gamma$1 & 11.45 (8.92) & R-L-M-90-$\Gamma$1 & 9.67 (8.41) & R-L-55-95-$\Gamma$1 & 9.31 (9.25) \\
          & D-N-60 & 11.57 (8.91) & R-L-55-90-$\Gamma$1 & 9.71 (9.04) & D-L-65 & 9.37 (8.29) \\
          & D-L-65 & 11.7 (8.37) & D-N-60 & 9.74 (10.16) & R-L-55-90-$\Gamma$1 & 9.4 (8.95) \\
          & D-L-60 & 11.71 (9.66) & R-L-M-95-$\Gamma$1 & 9.92 (8.46) & R-L-M-95-$\Gamma$1 & 9.44 (8.37) \\
          & R-L-55-95-$\Gamma$1 & 11.75 (9.14) & R-N-M-90-$\Gamma$1 & 10.26 (8.97) & R-N-M-90-$\Gamma$1 & 9.71 (8.59) \\
          & R-L-55-90-$\Gamma$1 & 11.77 (9.56) & R-L-50-90-$\Gamma$1 & 10.34 (8.58) & R-N-50-95-$\Gamma$1 & 10.31 (8.35) \\
          & R-L-M-90-$\Gamma$1 & 11.87 (9.09) & R-N-50-90-$\Gamma$1 & 10.34 (10.36) & R-N-M-95-$\Gamma$1 & 10.41 (9.08) \\
          & R-L-50-90-$\Gamma$1 & 11.93 (9.15) & R-N-M-95-$\Gamma$1 & 10.67 (9.71) & R-L-M-90-$\Gamma$1 & 10.44 (8.36) \\
          & R-L-50-95-$\Gamma$1 & 12.08 (9.45) & R-L-50-95-$\Gamma$1 & 10.96 (8.43) & R-L-50-95-$\Gamma$1 & 10.81 (9.15) \\
    \midrule
    1     & Model  & Mean (Std) & Model  & Mean (Std) & Model  & Mean (Std) \\
    \midrule
          & R-N-M-95-$\Gamma$1 & 6.9 (7.76) & R-N-50-90-$\Gamma$1 & 8.58 (10.17) & D-N-80 & 7.95 (8.36) \\
          & R-L-M-95-$\Gamma$1 & 7.3 (8.27) & R-N-50-95-$\Gamma$2 & 8.74 (9.86) & D-N-70 & 8.18 (8.55) \\
          & R-L-50-90-$\Gamma$2 & 7.42 (7.86) & D-L-55 & 8.8 (10.83) & R-N-55-90-$\Gamma$1 & 8.35 (8.13) \\
          & R-N-55-90-$\Gamma$1 & 7.87 (8.35) & D-N-70 & 8.86 (10.15) & R-N-M-90-$\Gamma$1 & 8.69 (8.57) \\
          & R-L-M-90-$\Gamma$1 & 7.89 (9.09) & D-L-70 & 9.19 (9.91) & R-N-50-95-$\Gamma$1 & 9.55 (8.27) \\
          & D-N-70 & 7.97 (7.54) & R-N-50-95-$\Gamma$1 & 9.9 (10.03) & D-L-65 & 9.6 (7.33) \\
          & R-N-M-90-$\Gamma$1 & 8.04 (8.32) & D-L-60 & 10.04 (10.6) & D-N-60 & 9.73 (8.8) \\
          & D-L-65 & 8.11 (7.92) & R-N-55-90-$\Gamma$1 & 10.23 (9.35) & D-N-65 & 9.74 (8.34) \\
          & D-L-60 & 8.4 (8.3) & D-L-50 & 10.24 (9.12) & D-L-50 & 10.42 (7.7) \\
          & R-L-55-90-$\Gamma$1 & 8.62 (8.22) & D-L-65 & 10.44 (9.05) & D-N-50 & 11.83 (7.95) \\
    \bottomrule
    \end{tabular}%
  \label{tab:small_instances_TW}%
\end{table}%

Table \ref{tab:small_instances_TW} demonstrates that the methods exhibiting superior cost-based robustness also perform remarkably well in terms of robustness against time window violations. Notably, all the methods listed in Table \ref{tab:small_instances_TW} yield almost equivalent levels of time window robustness, which is approximately 10 percent worse than the most conservative setting. It's worth noting that this outcome is quite interesting, given that the most conservative setting performs poorly in terms of cost due to its overly predicted demand values.

\subsection{Results for Large Instances}
For the Solomon instances involving 100 customers, we employ an ALNS heuristic to derive solutions for both CVRPTW and its robust version, RCVRPTW. The \texttt{ALNS} package already offers a comprehensive framework for CVRPTW, and we have adapted this framework for RCVRPTW as well.  

We use the two destroy operators $\Psi^- = \{random\_removal, string\_removal\}$ and two repair operators $\Psi^+ = \{greedy\_repair, regret\_repair\}$ in ALNS heuristic as described in Section \ref{sec:alns}. The initial weights of the operators are set to one, and following each iteration, these weights are updated by applying a decay rate of 0.8 to the existing weights. The new weights assigned to operators in the current iteration depend on the solutions they produce, whether they are the best so far, better, accepted or rejected, corresponding to new weights of 25, 5, 1 or 0, respectively. The initial temperature in the simulated annealing mechanism is tuned to ensure that there is a 0.05 chance of selecting a solution up to 50\% percent worse than the initial solution and the temperature is updated at each interaction such that it reaches 1 in 5,000 iterations. For each instance, run ALNS for 60 seconds.

%We adopt two destroy operators, $\Psi^- = \{random\_removal, string\_removal\}$. The $random\_removal$ operator removes a predefined number of customers from the routes of a feasible solution randomly. The $string\_removal$ operator removes partial routes around a randomly selected customer. As for repair operators, we utilize $\Psi^+ = \{greedy\_repair, regret\_repair\}$. The $greedy\_repair$ approach inserts unassigned customers into existing routes in a greedy manner. The $regret\_repair$ operator first creates a list of unassigned customers according to their decreasing regrets. Regret is determined by calculating the cost difference between assigning a customer to its best position and its second-best position. The $regret\_repair$ operator then inserts these unassigned customers from the list into the solution in a greedy manner, following the order of the list.

%The accept function in our heuristic is implemented using a simulated annealing mechanism, facilitating enhanced solution exploration at the beginning of the optimization process. The initial weights of the operators are set to 1, and following each iteration, these weights are updated by applying a decay rate to the existing and new weights. The new weights allocated to operators in the current iteration depend on the solutions generated by these operators, whether they are the best so far, better, accepted, or rejected solutions. ALNS is stopped after running and the best solution found in terms of initial cost is recorded.  

\begin{figure}[htbp]
    \centering
   \begin{subfigure}{0.332\textwidth}
        \centering
        \includegraphics[width=\linewidth]{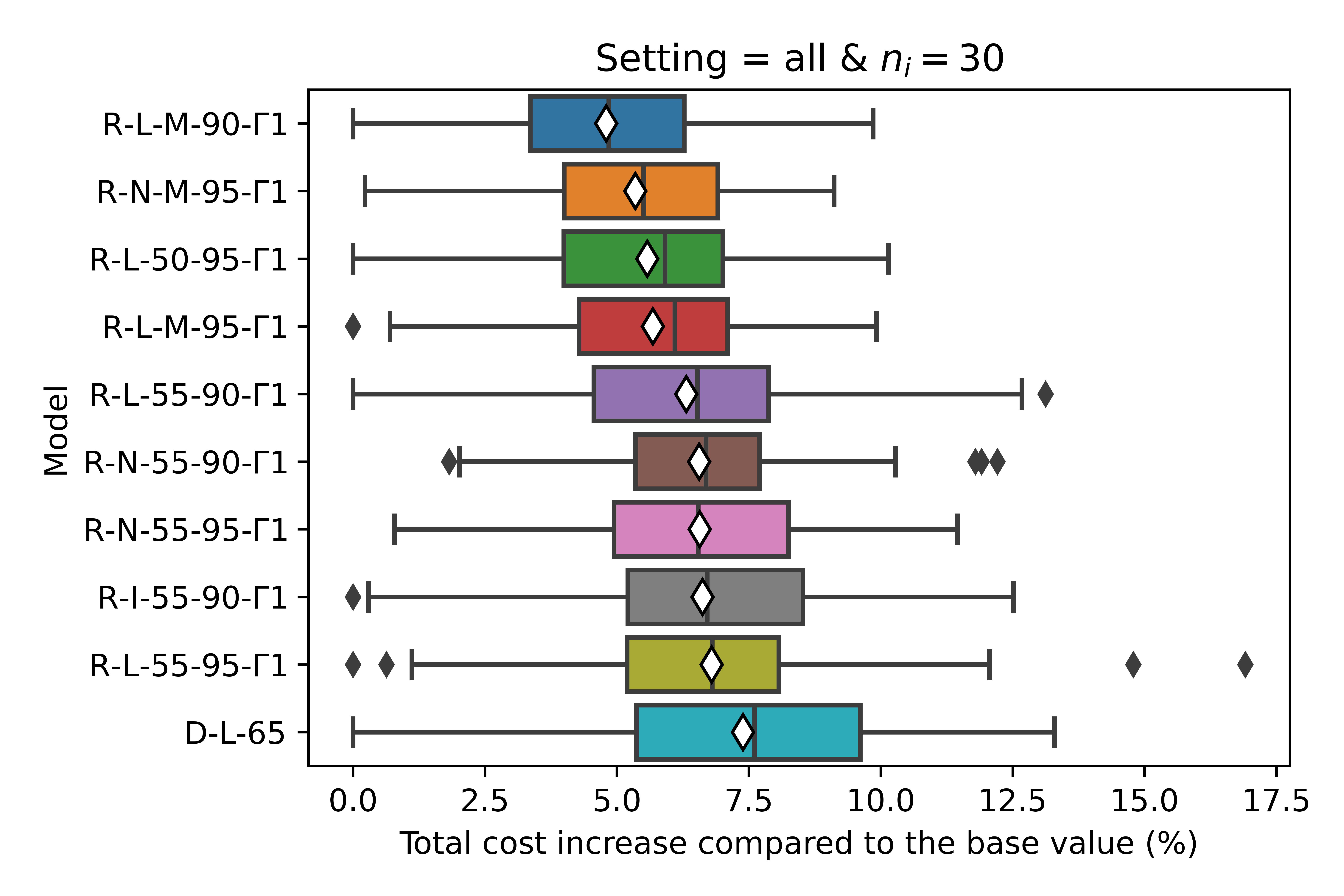}
    \end{subfigure} \hspace{-0.2cm} \begin{subfigure}{0.332\textwidth}
        \centering
        \includegraphics[width=\linewidth]{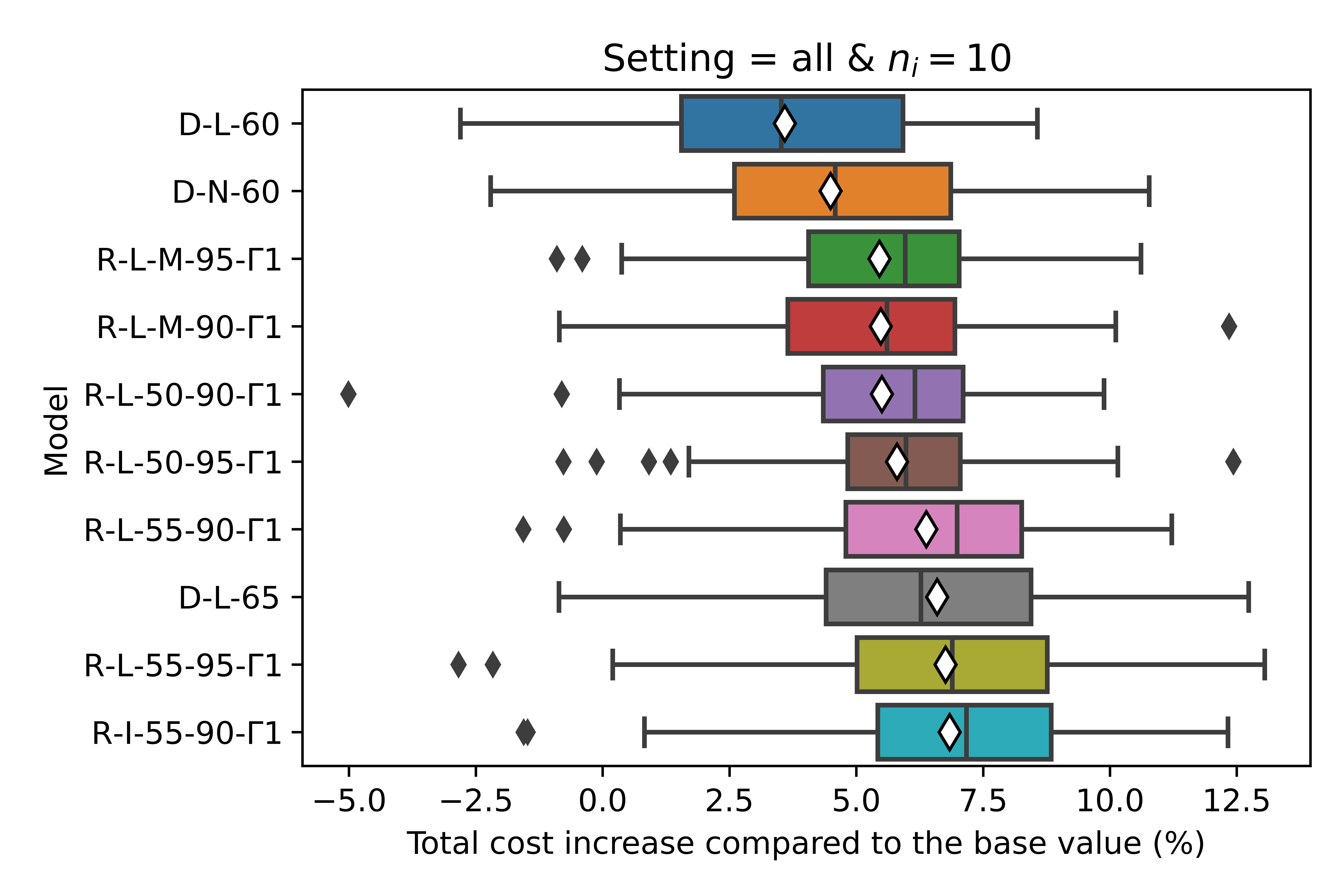}
   \end{subfigure} \hspace{-0.2cm}  \begin{subfigure}{0.332\textwidth}
        \centering
        \includegraphics[width=\linewidth]{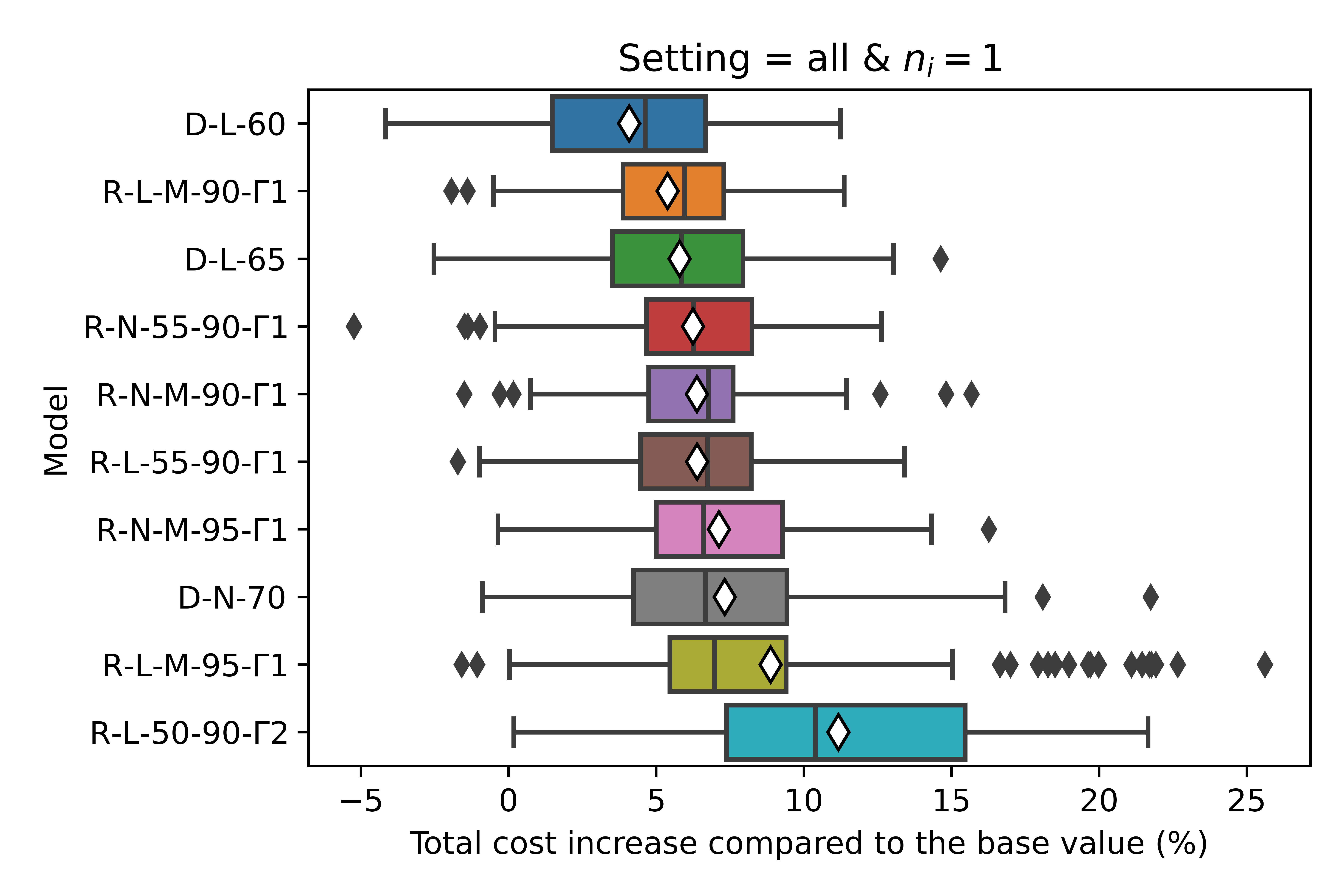}
   \end{subfigure}

    \vspace{0.5cm} % Add some vertical space between rows
    
    \begin{subfigure}{0.332\textwidth}
        \centering
        \includegraphics[width=\linewidth]{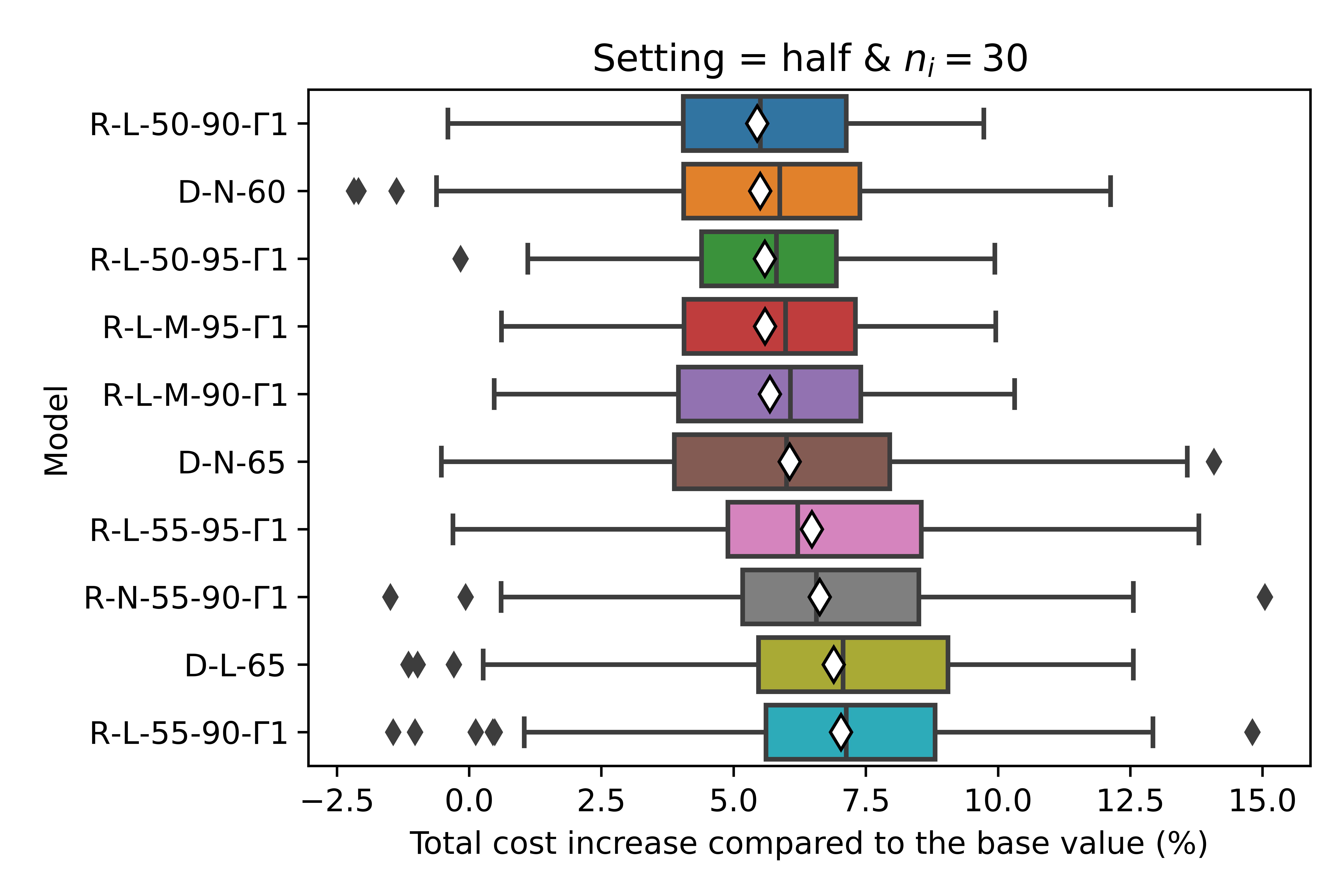}
           
    \end{subfigure}  \hspace{-0.2cm} \begin{subfigure}{0.332\textwidth}
        \centering
        \includegraphics[width=\linewidth]{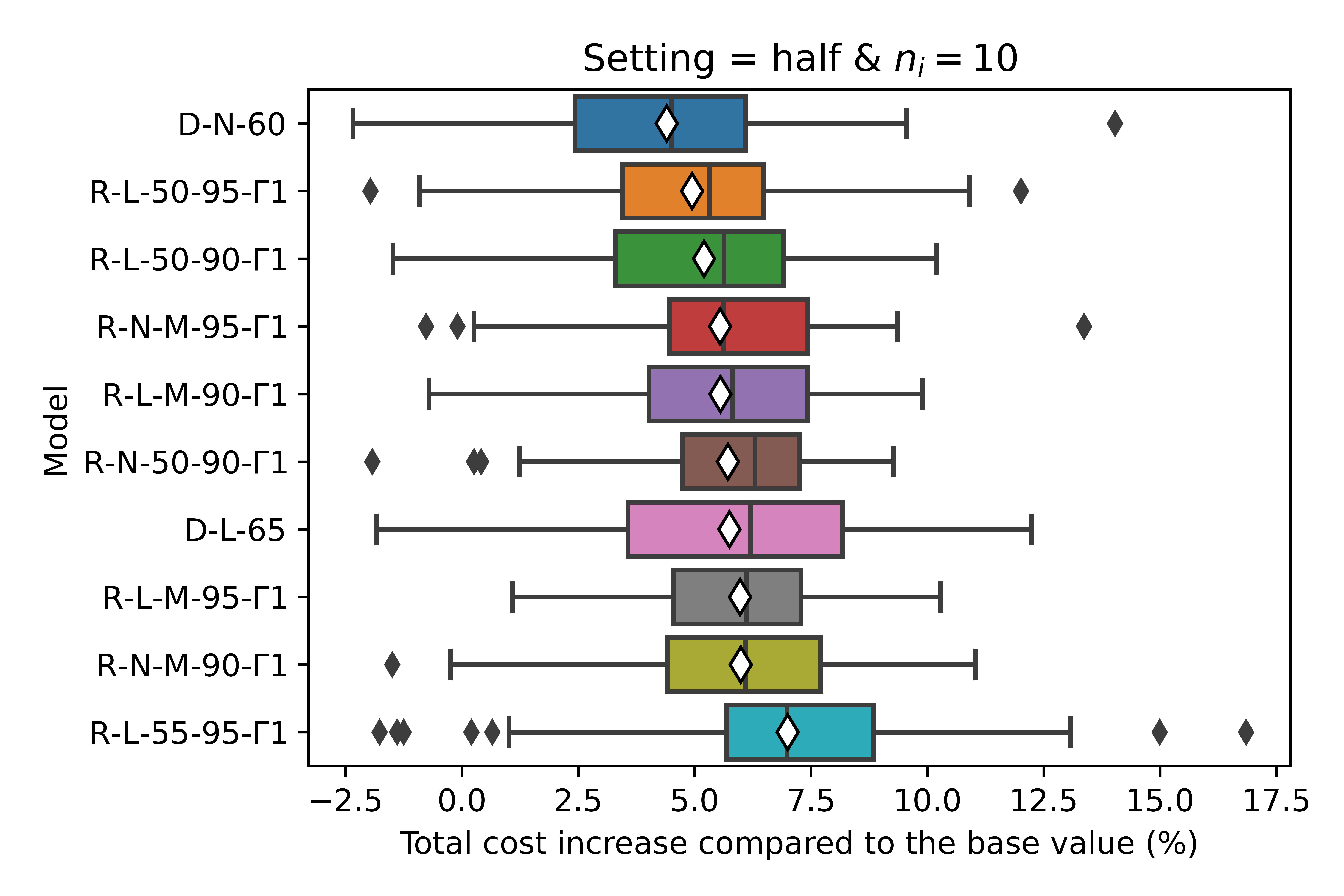}
         
    \end{subfigure} \hspace{-0.2cm}  \begin{subfigure}{0.332\textwidth}
        \centering
        \includegraphics[width=\linewidth]{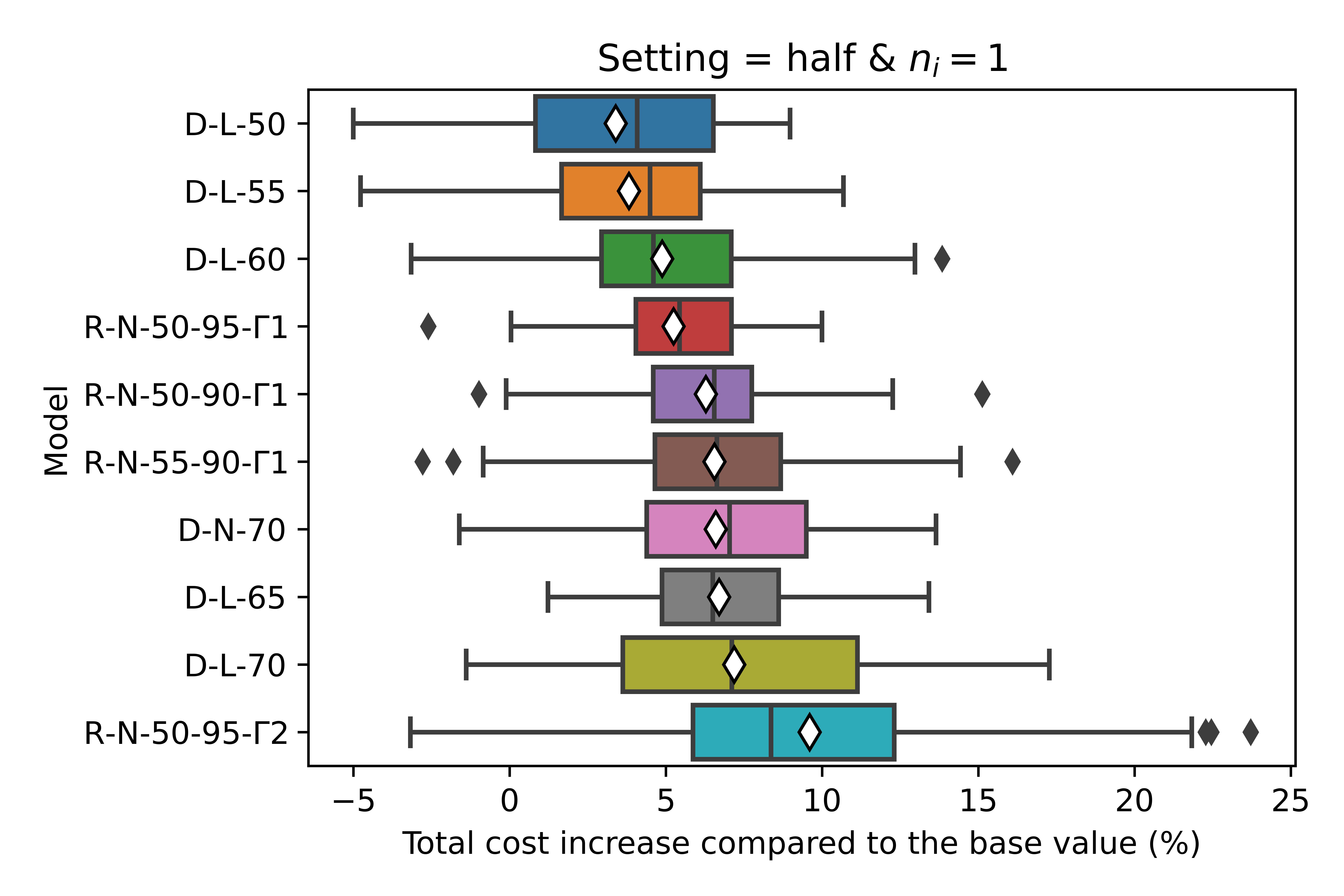}
    \end{subfigure}
    
    \vspace{0.5cm} % Add some vertical space between rows
    
    \begin{subfigure}{0.332\textwidth}
        \centering
        \includegraphics[width=\linewidth]{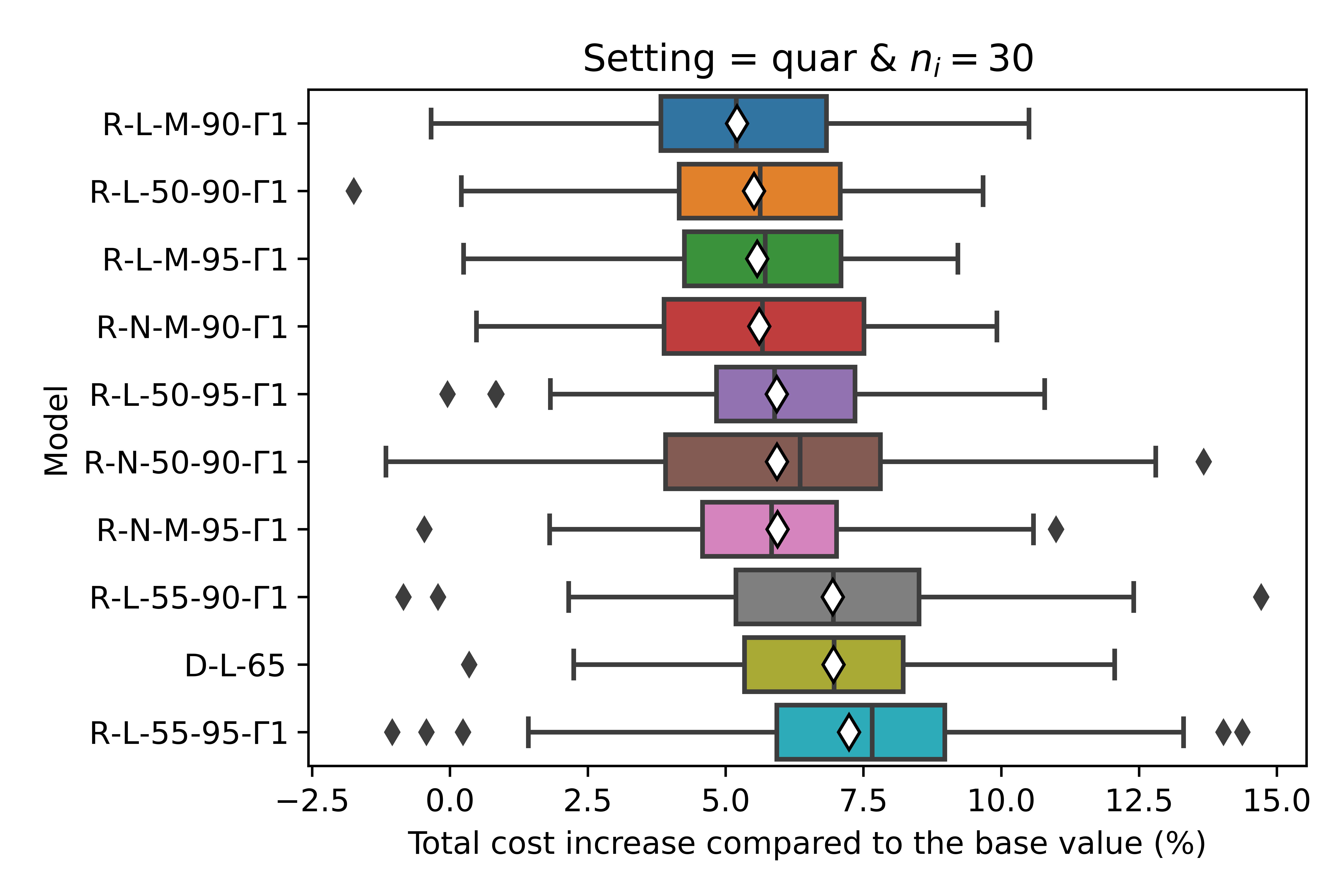}
          
    \end{subfigure} \hspace{-0.2cm} \begin{subfigure}{0.332\textwidth}
        \centering
        \includegraphics[width=\linewidth]{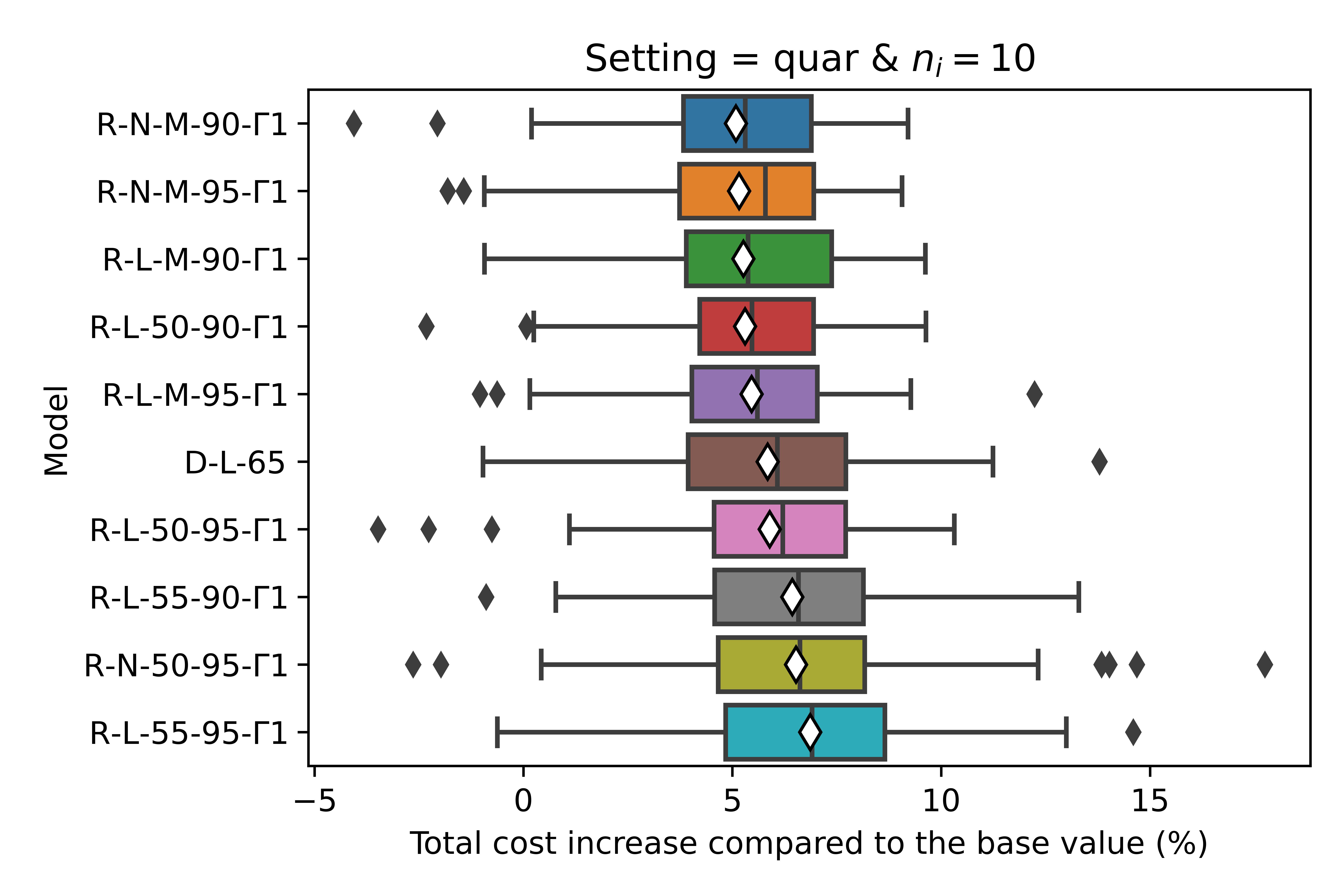}
         
    \end{subfigure}  \hspace{-0.2cm} \begin{subfigure}{0.332\textwidth}
        \centering
        \includegraphics[width=\linewidth]{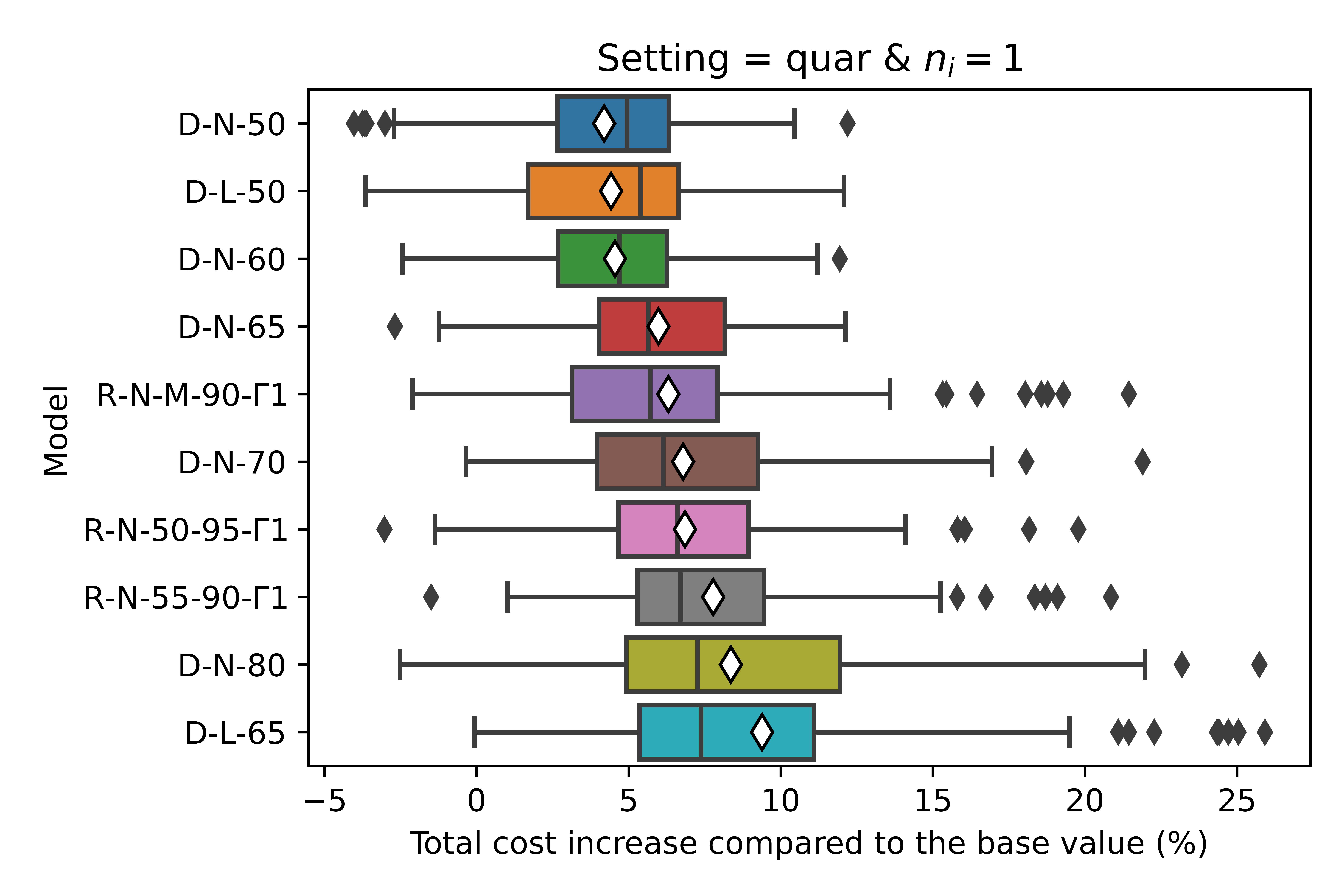}
          
    \end{subfigure}
    \caption{Performances of models for different settings instances with 100 customers. White and black diamonds represent the average values and outliers, respectively.}
    \label{fig:results_large_instances}
\end{figure}

Figure \ref{fig:results_large_instances} illustrates the simulation outcomes in terms of normalized cost, similar to the small instances. However, an intriguing observation emerges regarding the performance of deterministic and robust optimization models in exact and heuristic solution settings. Specifically, in richer data contexts, the relative performance of deterministic models deteriorates, whereas, in data-scarce scenarios, their relative performance improves.  

In the \textit{all} setting with $n_i = 30$, for example, the method D-L-65 ranks as the third best in terms of average performance in small problems. However, in the larger instances, it drops to the 10th best position. On the contrary, the top four methods are all deterministic models for large instances, whereas this was not the case for small instances. This result underscores the importance of utilizing complex prediction and optimization models in data-rich environments, as their advantages become less prominent when data scarcity is prevalent. 

Finally, we present the fractions of customers with time window violations for the larger instances in Table \ref{tab:large_instances_TW}. The table consistently demonstrates that the deterministic models generally result in fewer time window violations compared to the robust models. Although the difference is not substantial, this finding underscores another advantage of employing deterministic models with tailored demand quantile predictions to attain robust solutions, even though the model itself doesn't represent robustness explicitly.

% Table generated by Excel2LaTeX from sheet 'Sheet2'
\begin{table}[htbp]
  \centering
  \caption{Average percentage increase in the customers with time window violations compared to the most conservative model D-$\delta$-95 for each setting and number of observations ($n_i$) for large instances.}
    \begin{tabular}{c|cc|cc|cc}
    \multicolumn{1}{c}{} & \multicolumn{6}{c}{Setting} \\
\cmidrule{2-7}    $n_i$   & \multicolumn{2}{c}{\textit{all}} & \multicolumn{2}{c}{\textit{half}} & \multicolumn{2}{c}{\textit{quar}} \\
    \midrule
    \midrule
    30    & Model  & Mean (Std) & Model  & Mean (Std) & Model  & Mean (Std) \\
    \midrule
          & D-L-65 & 50.05 (23.84) & D-N-60 & 49.88 (24.4) & D-L-65 & 49.93 (24.14) \\
          & R-I-55-90-$\Gamma$1 & 51.15 (25.35) & D-N-65 & 49.95 (24.13) & R-L-55-90-$\Gamma$1 & 50.92 (25.91) \\
          & R-L-M-90-$\Gamma$1 & 51.33 (26.41) & D-L-65 & 49.99 (23.9) & R-L-55-95-$\Gamma$1 & 51.23 (25.04) \\
          & R-L-55-95-$\Gamma$1 & 51.41 (25.24) & R-N-55-90-$\Gamma$1 & 50.64 (24.77) & R-L-M-90-$\Gamma$1 & 51.64 (26.48) \\
          & R-L-55-90-$\Gamma$1 & 51.47 (25.72) & R-L-55-90-$\Gamma$1 & 50.81 (25.25) & R-N-50-90-$\Gamma$1 & 51.83 (26.41) \\
          & R-N-55-90-$\Gamma$1 & 51.75 (26.0) & R-L-55-95-$\Gamma$1 & 50.84 (25.23) & R-N-M-95-$\Gamma$1 & 51.96 (26.36) \\
          & R-L-50-95-$\Gamma$1 & 51.88 (26.46) & R-L-M-95-$\Gamma$1 & 51.74 (26.56) & R-L-50-95-$\Gamma$1 & 52.01 (26.49) \\
          & R-N-M-95-$\Gamma$1 & 51.91 (26.34) & R-L-50-90-$\Gamma$1 & 51.76 (26.65) & R-L-50-90-$\Gamma$1 & 52.2 (26.96) \\
          & R-N-55-95-$\Gamma$1 & 52.12 (26.24) & R-L-50-95-$\Gamma$1 & 52.23 (26.15) & R-N-M-90-$\Gamma$1 & 52.34 (26.67) \\
          & R-L-M-95-$\Gamma$1 & 52.17 (26.52) & R-L-M-90-$\Gamma$1 & 52.49 (26.53) & R-L-M-95-$\Gamma$1 & 52.44 (26.36) \\
    \midrule
    10    & Model  & Mean (Std) & Model  & Mean (Std) & Model  & Mean (Std) \\
    \midrule
          & D-L-60 & 49.98 (24.97) & D-N-60 & 49.99 (24.57) & D-L-65 & 49.81 (23.75) \\
          & D-L-65 & 50.18 (23.96) & D-L-65 & 50.12 (24.0) & R-L-55-90-$\Gamma$1 & 51.22 (25.56) \\
          & D-N-60 & 50.69 (25.25) & R-L-55-95-$\Gamma$1 & 51.47 (25.19) & R-N-50-95-$\Gamma$1 & 51.53 (25.29) \\
          & R-I-55-90-$\Gamma$1 & 51.04 (25.15) & R-N-M-90-$\Gamma$1 & 51.74 (26.59) & R-L-55-95-$\Gamma$1 & 51.64 (25.4) \\
          & R-L-55-95-$\Gamma$1 & 51.2 (25.65) & R-N-M-95-$\Gamma$1 & 51.85 (26.78) & R-L-M-90-$\Gamma$1 & 51.82 (26.34) \\
          & R-L-55-90-$\Gamma$1 & 51.32 (25.26) & R-L-50-95-$\Gamma$1 & 51.92 (26.37) & R-N-M-90-$\Gamma$1 & 51.84 (26.89) \\
          & R-L-50-90-$\Gamma$1 & 52.02 (26.59) & R-L-50-90-$\Gamma$1 & 51.95 (26.56) & R-N-M-95-$\Gamma$1 & 51.85 (26.66) \\
          & R-L-50-95-$\Gamma$1 & 52.09 (26.53) & R-L-M-90-$\Gamma$1 & 52.12 (26.23) & R-L-50-90-$\Gamma$1 & 51.98 (26.51) \\
          & R-L-M-95-$\Gamma$1 & 52.16 (26.72) & R-L-M-95-$\Gamma$1 & 52.21 (26.29) & R-L-50-95-$\Gamma$1 & 52.13 (26.4) \\
          & R-L-M-90-$\Gamma$1 & 52.18 (26.41) & R-N-50-90-$\Gamma$1 & 52.56 (26.61) & R-L-M-95-$\Gamma$1 & 52.31 (26.4) \\
    \midrule
    1     & Model  & Mean (Std) & Model  & Mean (Std) & Model  & Mean (Std) \\
    \midrule
          & D-N-70 & 49.47 (22.84) & D-L-70 & 49.0 (22.47) & D-N-70 & 49.18 (22.84) \\
          & D-L-65 & 49.62 (23.91) & D-N-70 & 49.8 (23.6) & D-N-80 & 49.57 (22.58) \\
          & D-L-60 & 50.01 (24.65) & D-L-60 & 50.11 (24.72) & D-N-60 & 50.16 (24.66) \\
          & R-L-55-90-$\Gamma$1 & 50.51 (24.89) & D-L-65 & 50.35 (24.44) & D-L-50 & 50.93 (27.13) \\
          & R-N-55-90-$\Gamma$1 & 50.65 (24.74) & D-L-55 & 50.44 (25.82) & D-N-65 & 50.99 (25.42) \\
          & R-N-M-95-$\Gamma$1 & 51.18 (25.59) & R-N-55-90-$\Gamma$1 & 51.16 (25.17) & D-L-65 & 51.07 (23.43) \\
          & R-N-M-90-$\Gamma$1 & 51.52 (25.62) & D-L-50 & 51.43 (27.28) & R-N-M-90-$\Gamma$1 & 51.12 (25.21) \\
          & R-L-M-90-$\Gamma$1 & 51.59 (26.46) & R-N-50-90-$\Gamma$1 & 51.86 (26.05) & R-N-50-95-$\Gamma$1 & 51.13 (25.26) \\
          & R-L-M-95-$\Gamma$1 & 51.81 (25.82) & R-N-50-95-$\Gamma$2 & 52.05 (25.38) & D-N-50 & 51.32 (26.71) \\
          & R-L-50-90-$\Gamma$2 & 52.2 (24.79) & R-N-50-95-$\Gamma$1 & 52.11 (26.06) & R-N-55-90-$\Gamma$1 & 51.62 (25.22) \\
    \bottomrule
    \end{tabular}%
  \label{tab:large_instances_TW}%
\end{table}%
 
\section{Conclusion} \label{sec:conclusion}
This paper investigates the utilization of contextual information to address the RCVRPTW, presenting its practical relevance when facing scenarios with limited demand history for certain customers. More specifically, we investigated the use of quantile prediction, both for making single predictions for use in deterministic optimization as well as for predicting the uncertainty set for robust optimization.  Our results also offer insights into robust decision-making, particularly for new customers lacking historical data. The approach chosen can vary depending on data availability. In data-rich contexts, multiple quantile predictions and robust optimization models are advantageous. Conversely, for data-scarce situations, single predictions and deterministic optimization models maintain strong performance in terms of cost and time window violations.

Predicting quantile demands offers an added advantage by enabling robust decision-making with the flexibility to adjust the level of robustness. Unlike mean demand predictions, quantile predictions allow the decision-maker to fine-tune the degree of robustness. These predicted quantiles can be integrated into both exact solution methods and heuristic approaches.

The deterministic model is the industry standard and it is widely used by practitioners. Our results offer valuable insights for practitioners, showcasing the use of deterministic models, which are often more prevalent than their robust counterparts. These deterministic models can provide commendable performance in terms of robustness, making them a viable option for non-experts in the field who seek user-friendly solutions without compromising on robustness.

The present study addresses prediction and optimization separately, even though predictions are made with the awareness that the predicted parameters will be used in an optimization model. An avenue for expanding this research involves integrating the prediction and decision-making processes in a single end-to-end pipeline (see, for example, \cite{mandi_decision-focused_2023} for a survey of the state-of-the-art end-to-end methods for prediction and optimization). While such a predict-then-optimize framework may pose computational challenges, it may potentially improve routing decisions since training loss is directly defined in terms of the decision quality of the downstream optimization problem. Furthermore, in an end-to-end pipeline, the decision-maker is not required to predefine a fixed quantile before solving the optimization problem to achieve the best performance.

\section*{Acknowledgements}
This research received funding from the
European Research Council (ERC) under the European Union’s Horizon 2020 research and innovation program (Grant No. 101002802, CHAT-Opt).
 
\clearpage

\bibliographystyle{apalike} 
{\footnotesize\bibliography{references,references1}}

\end{document}